\documentclass[final,onefignum,onetabnum]{siamart250211}

\ifpdf
\hypersetup{
	pdfauthor={Cyprien Tamekue, ShiNung Ching},
	pdfsubject={Paper manuscrit},
	pdftitle={Control analysis and synthesis for general control-affine systems},
	pdfkeywords={Control-affine systems, global controllability, control synthesis, Kalman’s type conditions, Gramian controllability matrix}
}
\fi

\graphicspath{{./figs/}}

\usepackage{bbm}
\usepackage{graphicx, amssymb,stmaryrd}
\usepackage{xcolor}
\usepackage{mathrsfs}
\usepackage{hyperref}
\usepackage{marvosym}
\usepackage{comment}
\usepackage{amsfonts, amssymb, amsmath}
\usepackage{mathtools}

\usepackage{tikz,tkz-tab} 
\usepackage{moreverb}

\usepackage{csquotes}

\numberwithin{equation}{section}
\newcommand{\bs}{\symbol{92}}

\newcommand*\samethanks[1][\value{footnote}]{\footnotemark[#1]}

\newsiamremark{remark}{Remark}
\newsiamremark{hypothesis}{Hypothesis}
\newsiamremark{example}{Example}
\newsiamremark{notation}{Notation}
\newsiamremark{assumption}{Assumption}
\crefname{hypothesis}{Hypothesis}{Hypotheses}

\def\idty{{\operatorname{Id}}}

\newcommand{\Ran}{\operatorname{Ran}}

\def\cF{{\mathcal F}}
\def\cZ{{\mathcal Z}}
\def\N{{\mathbb N}}

\def\R{{\mathbb R}}

\def\cA{{\mathcal A}}    \def\cS{{\mathcal S}} \def\cB{{\mathcal B}}   \def\cN{{\mathcal N}}      \def\cU{{\mathcal U}}      \def\cE{{\mathcal E}}     \def\cF{{\mathcal F}}  \def\cL{{\mathcal L}} \def\cR{{\mathcal R}} \def\cX{{\mathcal X}} \def\cY{{\mathcal Y}}  \def\cZ{{\mathcal Z}}

\headers{Control analysis and synthesis for general control-affine systems}{C. Tamekue and S. Ching}
\title{Control analysis and synthesis for general control-affine systems\thanks{\funding{This work is partially supported by grant R21MH132240 from the US National Institutes of Health to SC.}}}
\author{Cyprien Tamekue\thanks{Department of Electrical and Systems Engineering, Washington University in St. Louis, St. Louis, 63130, MO, USA. (\url{cyprien@wustl.edu}, \url{shinung@wustl.edu}). The first author thanks his former advisor, Yacine Chitour, for many relevant ideas and advice.}  \and ShiNung Ching\samethanks}

\begin{document}

\maketitle

\begin{abstract}
We study controllability and constructive synthesis for control–affine systems. We introduce trajectory–dependent Gramian maps that extend the linear time–varying Gramian and yield explicit fixed–point synthesis maps. On feasible coercivity classes (uniform eigenvalue lower bounds), the Gramian map is Lipschitz, and under a comparison estimate criterion, synthesis iterates exhibit decay, and the Banach fixed-point theorem gives a unique fixed point that steers the system and satisfies an energy identity. When, in addition, an orthogonality condition holds, this fixed point coincides with the unique global minimum–energy control on the feasible set; if the coercivity bound holds uniformly for all bounded controls, the same conclusion holds on the full bounded–control space. We provide structural conditions on the input matrix that ensure the nonemptiness of the feasible class (and, in fully actuated regimes, equality with the full space) and sufficient conditions for underactuated systems via bounded–amplitude reference controls. Case studies on Hopfield network dynamics illustrate refined estimates that enlarge reachable targets. A trajectory–freezing and compactness step extends the synthesis to general nonlinear control-affine systems. The results yield verifiable controllability criteria with explicit, numerically implementable controllers.
\end{abstract}

\begin{keywords}
Control-affine systems; nonlinear controllability; fixed-point methods; controllability Gramians; reachability; minimal-energy control.
\end{keywords}

\begin{MSCcodes}
93B05, 93C10, 93C15, 49K15, 49N35.
\end{MSCcodes}

\section{Introduction}\label{s:introduction}
We study control analysis and implementable synthesis for the control–affine class
\begin{equation}\label{eq:nonlinear control system general}
    \dot{x}(t)=A(t,x(t))\,N_t(x(t))+B(t,x(t))\,u(t),\qquad x(t_0)=x^0,\quad t\in[t_0,T],
\end{equation}
where at time $t$, the system's state is $x(t) = (x_1(t), \ldots, x_d(t))^\top \in \R^d$, $x^0 = (x_1^0, \ldots, x_d^0)^\top \in \R^d$ is the initial state, $u(t) = (u_1(t), \ldots, u_k(t))^\top \in \R^k$ is the control input, $N_t$ is a nonautonomous, vector field, $A(t,x)$ and $B(t,x)$ are time- and state-dependent matrices. Models of this form cover, among others, recurrent neural dynamics~\cite{hopfield1984neurons,sontag1998remarks} and standard engineering systems~\cite{cannon2003dynamics}, and hence their study remains important. 

Systems of the form~\eqref{eq:nonlinear control system general} have been extensively studied in control theory, with a wide range of mathematical techniques developed over the past decades.  
Classical approaches include linearization-based analysis \cite{sontag2013mathematical,coron2007control}, homotopy continuation \cite{chitour2006continuation,ji2023global,zhengping2024regularized,sussmann1992new}, and the return method \cite{coron2002return,coron2007control}, which remains one of the most powerful constructive tools for global controllability.  
Geometric and algebraic methods based on Lie algebras provide structural controllability criteria \cite{hirschorn1976global,lukes1972global,sarychev2006lie,sussmann1983lie,sussmann1987general,coron2007control}, while topological arguments offer complementary perspectives \cite{furi1985topological,sontag1998remarks}.  
Other analytic techniques include power series expansions \cite{coron2007control} and fixed-point approaches \cite{zuazua1991exact,zuazua1993exact,coron2007control,lukes1972global}. These contributions have led to a rich theory of local and global controllability; yet, implementable minimum-energy controls remain scarce beyond the linear time-invariant (LTI) and linear time-varying (LTV) settings~\cite{kalman1963controllability}, where the controllability Gramian yields closed-form minimum-energy controls.

Beyond the LTI/LTV setting, explicit Gramian–based formulas for minimum–energy controls are rarely available, and many nonlinear controllability results are qualitative or rely on continuation along a state path. Our stance is to recover ``Gramian calculus'' for system~\eqref{eq:nonlinear control system general} by introducing \emph{trajectory–dependent} controllability Gramian maps and casting synthesis as a fixed–point problem on an \emph{feasible coercivity class}. This yields (i) a constructive control law via an implicit Gramian formula; (ii) a built–in \emph{energy certificate}; and (iii) a convergence mechanism (via explicit Volterra rank-one comparison estimates for the iterates of the synthesis map), which is directly implementable. In contrast, the homotopy continuation method (HCM)~\cite{chitour2006continuation,ji2023global,sussmann1992new} advances by integrating the path lifting equation and therefore requires surjectivity of the differential of the endpoint map along the entire path and typically offers no intrinsic energy certificate. Note, however, that a regularized continuation method has been proposed in~\cite{schmoderer2022study} to deal with the case where this differential is not surjective. Still, the approach requires some conditions for the solvability: the drift dynamics must be well-posed over $[t_0, T]$, and the family
of regularized path lifting equations must converge to a solution of the original problem. See also~\cite{zhengping2024regularized,schmoderer2024regularised}. Our framework trades this global rank condition for verifiable, local requirements: the coercivity of a trajectory–dependent Gramian, a self–mapping inequality on the feasible class, and an explicit iterate-contractivity condition.
In practice, this matters when complete controllability is not known a priori: the method still furnishes steering controls (and their cost) whenever the feasible class is nonempty, and it delivers a quantitative subset of the reachable set. When HCM does apply, the two approaches are complementary---the nonlinear Gramian recovers the classical solution in the linear case and can serve as a preconditioner or initializer for continuation---while in regimes where the rank fluctuates, the Banach fixed–point route remains robust by construction.

Our main contributions are the following: 
\textit{(i) Baseline synthesis on $\dot x=N_t(x)+B(t,x)u$:} We work on a feasible coercivity class $\cF(C)$ (uniform lower bound on the trajectory–dependent Gramian). On $\cF(C)$, the synthesis operator is well defined, and---on a self–mapped ball---a Volterra rank-one comparison estimate gives an iterate-contractivity criterion, which yields a unique fixed point with the Picard convergence and an energy certificate. \textit{(ii) From relative to global minimality:} If an orthogonality condition holds at the fixed point, then this fixed point is the \emph{unique global} $L^2$–minimizer over the feasible set. \textit{(iii) Structural guarantees:} An integral nondegeneracy of the input matrix $B$ implies uniform coercivity ($\cF(C)=L^\infty$). More generally, a reference control $u^\flat$ with invertible Gramian and amplitude small relative to a prescribed model-driven constant ensures the self–map property for explicit target radii; model structure (e.g., Hopfield networks) sharpens constants and enlarges admissible targets. \textit{(iv) Reachability:} Self–mapping on $\cF(C)$  together with the iterate-contractivity condition provides constructive steering with quantitative energy bounds; if $\cF(C)=L^\infty$ one recovers complete controllability, otherwise one obtains a computable subset of the reachable set. \textit{(v) Extension to the general system \eqref{eq:nonlinear control system general}:} Freezing along an auxiliary trajectory reduces \eqref{eq:nonlinear control system general} to a family of baseline problems with uniform Lipschitz control of the Gramian map; Schaefer’s theorem then yields steering controls for admissible targets. For $N_t(x)=x$, the construction recovers \cite[Theorem.~3.40]{coron2007control}.

The remainder of the paper is structured as follows. In Notation~\ref{not:general} and Assumptions~\ref{ass:general assumption nonautonomous}, we introduce the notation and assumptions used throughout. Section~\ref{s:state-dependent-input} develops the synthesis framework for the baseline system $\dot x=N_t(x)+B(t,x)u$: we introduce the nonlinear Gramian representations, construct the synthesis maps, define the feasible coercivity class, and use an iterate-contractivity condition to obtain steering controls. Structurally sufficient conditions are then established, and refined estimates are illustrated on Hopfield-type networks.  
Section~\ref{s:sufficient conditions for nonautonomous quadratic nonlinear control systems with a pertubation} extends the framework to the general control-affine system~\eqref{eq:nonlinear control system general} by freezing along auxiliary trajectories and invoking a compactness argument to recover reachability for admissible targets. Section~\ref{s:conclusion} summarizes the results and discusses perspectives for future work, and main proofs are collected in Section~\ref{app:proofs of some results presented in the main text}.

\begin{notation}\label{not:general}   
In the following,  we denote by $\R^{d\times k}$ the set of $d\times k$ matrices with real coeﬃcients.  $\R^d$ denotes the vector of $d$ real columns in dimension and $\cL(\R^k, \R^d)$ the space of linear maps from $\R^k$ to $\R^d$ that we identify, in the usual way, with $\R^{d\times k}$ . For all $(x,y)\in\R^d\times\R^d$, we denote by $|x|$  the Euclidean norm of $x$ and by $\langle x,y\rangle$ the scalar product of $x$ and $y$. For a symmetric positive definite matrix $M\in\R^{d\times d}$, the norm of $x\in\R^d$ with respect to $M$ is defined by $|x|_M^2=x^\top Mx$. We denote the identity matrix of any size by $\idty$, and for every $M\in\R^{d\times k}$, $\|M\|=\sup\{|Mx|: x\in\R^k,\;|x|=1\}$ denotes the spectral norm of $M$. For a symmetric matrix $M\in\R^{d\times d}$, $\lambda_{\min}(M)$ and $\lambda_{\max}(M)$ denote respectively its smaller and largest eigenvalues. Finally, we use $\mathscr{L}(X, Y)$ to denote the space of linear bounded operators from two normed vector spaces $X$ and $Y$.
\end{notation}

\begin{assumption}\label{ass:general assumption nonautonomous} 
    Throughout the following, unless otherwise stated, the nonautonomous vector field $N_t:[t_0, T]\times\R^d\to\R^d$ satisfies the following assumptions
    \begin{enumerate}
        \item The map $t\mapsto N_t(x)$ belongs to $L^\infty(t_0, T)$ for every fixed $x\in\R^d$,
        \item The map $x\mapsto N_t(x)$ is $C^2$ for every fixed $t\in[t_0, T]$.
    \end{enumerate}
   Additionaly, for every $t\in[t_0, T]$, $N_t$ is globally $\Lambda_1$-Lipschitz function for some $\Lambda_1>0$, viz.
   \begin{equation}\label{eq:EGectral norm of DN_t}
    \|DN_t(w)\|\le\Lambda_1\qquad\forall(t,w)\in[t_0, T]\times\R^d,
\end{equation}
where $DN_t(w)\in\R^{d\times d}$ is the differential of $N_t$ at any $w\in\R^d$.  Finally, we also assume that its second derivative $D^2N_t$ satisfy for some $\Lambda_2>0$
    \begin{equation}\label{eq:EGectral norm of D^2N_t}
        \|D^2N_t(x)\|\le\Lambda_2\qquad\forall (t, x)\in[t_0, T]\times\R^d.
    \end{equation}
\end{assumption}

\begin{remark}
   For every fixed $t\in[t_0, T]$, the map $x\mapsto N_t(x)$ does not need to be bounded. Throughout the manuscript, we let $\Delta t_0:=T-t_0$.
\end{remark}

\section{Prerequisites and flow representation}\label{s:prerequisites}
It follows from Assumption~\ref{ass:general assumption nonautonomous} that the nonautonomous flow $t\in[t_0, T]\mapsto\Phi_{t_0,t}$ associated with $N_t$ is globally Lipschitz 
and solves the following equation \cite[Chapter 2]{bressan2007introduction}
\begin{equation}\label{eq:nonautonomous flow}
    \frac{d}{d t}\Phi_{t_0,t}(x^0) = N_t(\Phi_{t_0,t}(x^0)),\qquad \Phi_{t_0,t_0}(x^0) = x^0\in\R^d.
\end{equation}
Moreover, $\{\Phi_{s,t}\mid (s,t)\in[t_0, T]^2\}$ forms a two-parameter family of diffeomorphisms which satisfies the following algebraic identities
\begin{align}\label{eq:algebraic identities}
    \Phi_{t,t} &= \idty\qquad\forall t\in[t_0, T],\\
     \Phi_{t_2,t_3}\circ\Phi_{t_1,t_2} &= \Phi_{t_1,t_3}\qquad\forall (t_1,t_2,t_3)\in[t_0, T]^3,\label{eq:algebraic identities 2}\\
     (\Phi_{t_1,t_2})^{-1}&=\Phi_{t_2,t_1}\qquad\forall (t_1,t_2)\in[t_0, T]^2.\label{eq:identity algebraic}
\end{align}
Furthermore, for fixed $(s,t)\in[t_0, T]^2$, the maps $\Phi_{s,t}$ and $\Phi_{t,s}$ 
are differentiable at every $x^0$ and $\Phi_{s,t}(x^0)$, respectively. Denoting these differential by $D\Phi_{s,t}(x^0)$, and $D\Phi_{t, s}(\Phi_{s,t}(x^0))$, 
they solve
\begin{equation}\label{eq:equation derivative of phi_t_0,t}
           \begin{cases}
            \displaystyle\frac{d}{d t} D\Phi_{s,t}(x^0)&=DN_t(\Phi_{s,t}(x^0))\,D\Phi_{s,t}(x^0),\qquad D\Phi_{s,s}(x^0)=\idty,\\
            \cr
             \displaystyle\frac{d}{dt} D\Phi_{t,s}(\Phi_{s,t}(x^0))&=-D\Phi_{t,s}(\Phi_{s,t}(x^0))\,DN_t(\Phi_{s,t}(x^0)),\qquad D\Phi_{s,s}(\Phi_{s,t}(x^0))=\idty.
        \end{cases}
    \end{equation}
In particular, the following holds for every $x^0\in\R^d$,
\begin{equation}\label{eq:inverse Dphi_t_0,t}
    [D\Phi_{s,t}(x^0)]^{-1} = D\Phi_{t,s}(\Phi_{s,t}(x^0))\qquad\forall (s,t)\in[t_0, T]^2.
\end{equation}

We collect useful estimates related to $\Phi_{s, t}$ and $\Phi_{t, s}$ in the following lemma. The proof uses the standard Gronwall lemma.

\begin{lemma}\label{lem:spectral norm of D_s_t_phi and D^2_s_t_phi}
    Let $\beta\in C^0([t_0, T];\R^d)$ and set $q:=\Lambda_2/\Lambda_1$. It holds for every $(s, t)\in[t_0, T]^2$, $s\le t$,
    \begin{equation}\label{eq:EGectral norm of D_s_t_phi}
    \begin{split}
         \|D\Phi_{s, t}(\beta(s))\|&\le e^{\Lambda_1(t-s)},\qquad  \|D\Phi_{t,s}(\beta(s))\|\le e^{\Lambda_1(t-s)},\\
          \|D\Phi_{s, t}(\beta(t))\|&\le e^{\Lambda_1(t-s)},\qquad  \|D\Phi_{t,s}(\beta(t))\|\le e^{\Lambda_1(t-s)},
    \end{split}
\end{equation}
\begin{equation}\label{eq:EGectral norm of D^2_s_t_phi}
\begin{split}
    \|D^2\Phi_{s,t}(\beta(s))\|&\le q(e^{2\Lambda_1(t-s)}-e^{\Lambda_1(t-s)}),\qquad \|D^2\Phi_{t,s}(\beta(s))\|\le q(e^{2\Lambda_1(t-s)}-e^{\Lambda_1(t-s)}),\\
    \|D^2\Phi_{s,t}(\beta(t))\|&\le q(e^{2\Lambda_1(t-s)}-e^{\Lambda_1(t-s)}),\qquad \|D^2\Phi_{t,s}(\beta(t))\|\le q(e^{2\Lambda_1(t-s)}-e^{\Lambda_1(t-s)}),
\end{split}
\end{equation}
where $D^2\Phi_{s,t}$ is the second derivative (third-order tensor) of $\Phi_{s,t}$. 
\end{lemma}

\section{Control–affine dynamics with state–dependent and time-varying input matrix}\label{s:state-dependent-input}
In this section, we consider the following control–affine system with state–dependent and time-varying input matrix
\begin{equation}\tag{\textrm{$\Sigma$}}\label{eq:state-dependent-input}
    \dot{x}(t) = N_t(x(t))+B(t,x(t))u(t),\quad x(t_0) = x^0,\quad t\in[t_0, T],
\end{equation}
which is a specific case of \eqref{eq:nonlinear control system general} where $A(t,x)=\idty$. We recall that $(t,x)\mapsto N_t(x)$ satisfies Assumption~\ref{ass:general assumption nonautonomous}. 
\begin{assumption}\label{ass:gneneral assumption on B}
   The input matrix $B:[t_0, T]\times\R^d\to\R^{d\times k}$ is an element of $L^\infty((t_0, T)\times\R^d; \R^{d\times k})$. Furthermore, the map $x\mapsto B(t,x)$ is $C^1$ and globally Lipschitz for every $t\in[t_0, T]$, i.e., there exists $L_B\ge 0$ such that
\begin{equation}
    \|D_xB(t,w)\|\le L_B\qquad\forall(t,w)\in[t_0, T]\times\R^d.
\end{equation}
Note that if $B(t,\cdot)=B(t)$, then $L_B = 0$.
\end{assumption}

\subsection{Representation of solutions}\label{ss:representation of solutions}

The existence and uniqueness of global solutions to \eqref{eq:state-dependent-input} follow as a straightforward application of the Cauchy-Lipschitz theorem, given Assumptions~\ref{ass:general assumption nonautonomous} and~\ref{ass:gneneral assumption on B}, see, for instance, \cite[Theorem 3.2.1]{bressan2007introduction}. However, the main objective of this paper is to synthesize a control function that solves the associated two-point boundary value problem. To this end, leveraging the regularity properties of the drift $N_t$, we consider solution representations of \eqref{eq:state-dependent-input} that facilitate control design. These solution representations fall within the general framework of chronological calculus introduced in \cite{agrachev1979exponential} and leveraged in \cite{kawski2011chronological,sarychev2006lie} for applications in geometric control theory. See also \cite[Chapter~6]{agrachev2019comprehensive} for a comprehensive presentation of this theory. In \cite{tamekue2024control},  solution representations of nonlinear control systems like \eqref{eq:state-dependent-input} were established when $N_t=N$, $B(t,x)=B$, and the control was constant; see also the proof of~\cite[Theorem~4.4]{tamekue2024mathematical}. In contrast, the present work extends these representations to the framework of \eqref{eq:state-dependent-input}. For completeness, we replace in \eqref{eq:state-dependent-input}, $B(t,x)\,u(t)$ with $b(t,x)\in L^\infty((t_0, T)\times\R^d; \R^d)$ where $x\mapsto\,b(t,x)$ is locally Lipschitz. The proof of the following theorem is given in Section~\ref{ss:proof of sol representation nonautonomous}.
\begin{theorem}\label{thm:sol representation nonautonomous}
For all $x^0\in\R^d$, and $b\in L^\infty((t_0, T)\times\R^d; \R^d)$ such that  $x\mapsto\,b(t,x)$ is locally Lipschitz, the system
$\dot{x}(t)=N_t(x(t))+b(t,x(t))$, $x(t_0)=x^0$ admits a unique absolutely continuous solution $x\in C^0([t_0,T];\R^d)$ that can be written for $i\in\{1,2\}$,
\begin{equation}\label{eq:forward-backward representation}
x(t)=\Phi_{\tau_i,t}\!\left(\Phi_{t_0,\tau_i}(x^0)+\int_{t_0}^{t} D\Phi_{s,\tau_i}\big(x(s)\big)\,b(s, x(s))\,ds\right),
\qquad t\in[t_0,T], \quad \tau_1=t_0,\;\tau_2=T.
\end{equation}
\end{theorem}

The representation corresponding to $\tau_1=t_0$ is termed \textit{forward representation}, while that corresponding to $\tau_2=T$ is termed \textit{backward representation}, consistent with their temporal directionality.

Observe that if $N_t \equiv A(t)$, where $A : [t_0, T] \to \mathbb{R}^{d \times d}$ is a matrix-valued function in $L^\infty((t_0, T); \mathbb{R}^{d \times d})$, then~\eqref{eq:forward-backward representation} reduces to the classical representation of solutions for linear time-varying (LTV) systems. 

\begin{corollary}\label{cor:solution in the linear case}
If $N_t\equiv A(t)\in L^\infty((t_0, T); \R^{d\times d})$ and $b(t,\cdot)\equiv b(t)\in L^\infty((t_0, T); \R^d)$, then \eqref{eq:forward-backward representation} recasts as
\begin{equation}\label{eq:linear nonautonomous solution representation}
    x(t) = R(t,t_0)\,x^0+\int_{t_0}^{t}R(t,s)\,b(s)\,ds
    \qquad\forall t\in[t_0, T],
\end{equation}
where $R(t,s)\in C^0([t_0,T]^2;\R^{d\times d})$ is the state-transition matrix of $\dot{y}(t) = A(t)\,y(t)$ satisfying $R(s,s)=\idty$.
\end{corollary}
\begin{proof}
 First, by \eqref{eq:nonautonomous flow}-\eqref{eq:identity algebraic} and the definition of $R(t,s)$, one has $\Phi_{s, t}(x)=R(t,s)\,x$ for all $x\in\R^d$. It follows that $D\Phi_{s, t}(x(t))=R(t,s)$ for all $(t,s)\in[t_0, T]^2$. Therefore, from~\eqref{eq:forward-backward representation}, one deduces 
\begin{equation}
    x(t) = \Phi_{\tau_i,t}\left(\Phi_{t_0,\tau_i}(x^0)+\int_{t_0}^{t}\!\!D\Phi_{s,\tau_i}(x(s))\,b(s) \, ds\right)=R(t,t_0)\,x^0+R(t,\tau_i)\int_{t_0}^{t}\!\!R(\tau_i,s)\,b(s) \, ds
\end{equation}
which yields~\eqref{eq:linear nonautonomous solution representation}.
\end{proof}


We now establish results that relate the state-transition matrix of the linearized form of control system~\eqref{eq:state-dependent-input} to derivatives of the flow of the vector field $N_t$ along the trajectory. It will be essential for our subsequent results. The proof is presented in Section~\ref{ss:state-transition matrix of the nonlinear control}.

\begin{lemma}\label{lem:state-transition matrix of the nonlinear control B=0}
    Let $(x^0, u)\in\R^d\times L^\infty((t_0, T); \R^k)$, and let $x_u\in C^0([t_0, T]; \R^d)$ be the solution of~\eqref{eq:state-dependent-input}.
    Then, the state-transition matrix $R_u(t,s)\in C^0([t_0,T]^2;\R^{d\times d})$ of the linearized equation 
    \begin{equation}\label{eq:linearized equation B=0}
       \dot{y}(t) = \big[DN_t(x_u(t))+D_xB(t,x_u(t))u(t)\big] y(t),\qquad t\in [t_0, T]
    \end{equation}
    satisfying $R_u(s,s)=\idty$, is given by $R_u(t,s)= P_{\tau_i,u}(t,s)$, $i\in\{1,2\}$ where $P_{\tau_i,u}(t,s)$ can be factorized as
    \begin{equation}\label{eq:state-transition matrix of the nonlinear control}
     P_{\tau_i,u}(t,s) := \big[D\Phi_{t,\tau_i}(x_u(t))\big]^{-1}\,M_{\tau_i,u}(t,s)\,D\Phi_{s,\tau_i}(x_u(s)),\qquad(t,s) \in [t_0, T]^2,\quad\tau_1=t_0,\,\tau_2=T.
    \end{equation}
    Here, $M_{\tau_i,u}(t,s)\in C^0([t_0,T]^2;\R^{d\times d})$ satisfying $M_{\tau_i,u}(s,s)=\idty$ is the state-transition matrix of 
     \begin{equation}\label{eq:state-transition matrix M_T,u}
            \dot{y}(t) = \cA_{\tau_i,u}(t)\,y(t),\qquad t\in[t_0, T]
        \end{equation}
         where $\cA_{\tau_i,u}(t):=\big[D^2\Phi_{t,\tau_i}(x_u(t))B(t,x_u(t))u(t)+D\Phi_{t,\tau_i}\big(x_u(t)\big)D_xB\big(t,x_u(t)\big)u(t)\big]\big[D\Phi_{t,\tau_i}(x_u(t))\big]^{-1}\in\R^{d\times d}$.
\end{lemma}

\begin{remark}\label{rmk:matrix  S_{t_0,x^0}(s,t)}
   When $B(t,x)\equiv B(t)\in L^\infty((t_0, T); \R^k)$ is state-independent, it is clear from the factorization~\eqref{eq:state-transition matrix of the nonlinear control} that the state-transition matrix of $\dot y(t)=DN_t(x_u(t))y(t)$ is different from $D\Phi_{s,t}(x_u(s))$ in general, and that there is equality if the control $u=0$ or if $N_t$ is linear. In the former case, $M_{\tau_i,0}(t,s)=\idty$, and we use
    \begin{equation}\label{eq:matrix  S_{t_0,x^0}(s,t)}
    S_{t_0, x^0}(t,s):=P_{t_0,0}(t,s) = D\Phi_{t_0,t}(x^0)\,D\Phi_{s,t_0}(\Phi_{t_0,s}(x^0)) = D\Phi_{s,t}(\Phi_{t_0, s}(x^0))\qquad\forall(s,t)\in[t_0, T]^2
\end{equation}
    as the representation of the state-transition matrix of the linearized equation $\dot{y}(t) = DN_t(\Phi_{t_0,t}(x^0))\,y(t)$.
\end{remark}

\subsection{Controllability results}\label{ss:controllability results nonautonomous}

Given $(x^0,x^1)\in\R^d\times\R^d$, our goal in this section is to identify sufficient conditions under which system~\eqref{eq:state-dependent-input} can be steered from the initial state $x^0$ to the target $x^1$ over the horizon $[t_0,T]$. Whenever these conditions are met, we provide a synthesis of two feasible controls realizing this transfer.

\subsubsection{Synthesis, analysis, and optimal control results}\label{sss:Control synthesis and analysis results} As announced, we present our control synthesis, analysis, and optimal control results in this section. We start with recalling the following
\begin{definition}
Given $x^0\in\R^d$, we say that $x^1\in\R^d$ is \emph{reachable on} $[t_0,T]$ from $x^0$ if there exists a control $u\in L^\infty((t_0,T);\R^k)$ such that the  solution $x_u(\cdot)$ of \eqref{eq:state-dependent-input} satisfies $x_u(T)=x^1$. The reachable set is denoted by
\begin{equation}\label{eq:reachable set nonautonomous}
    \cR_{\Sigma}(T, x^0) := \Big\{\, x_u(T)\;:\; x_u(\cdot)\ \text{solves \eqref{eq:state-dependent-input} with }x(t_0)=x^0\,\Big\}.
\end{equation}
If $\cR_{\Sigma}(T, x^0)=\R^d$ for every $x^0\in\R^d$, one says that \eqref{eq:state-dependent-input} is \textit{completely controllable on} $[t_0, T]$.
\end{definition}

For a fixed $x^0\in\R^d$ and control $u\in L^\infty((t_0,T);\R^k)$, let $x_u(\cdot):=x_{u,x^0}(\cdot)$ be the corresponding trajectory of~\eqref{eq:state-dependent-input}. We associate with $u$ and $x^0$ the matrices ($\tau_1=t_0$ and $\tau_2=T$)
\begin{equation}\label{eq:nonautonomous N_i}
    \cN_{\tau_i}(u):=\cN_{\tau_i}(u, x^0) = \int_{t_0}^{T}D\Phi_{t,\tau_i}(x_u(t))\,B(t,x_u(t))\,B(t,x_u(t))^\top\,D\Phi_{t,\tau_i}(x_u(t))^\top\,dt,\quad\,i\in\{1,2\}.
\end{equation}

\begin{remark}\label{rmk:on the matrices N_1 and N_2}
(1) If $N_t\equiv A(t)\in L^\infty((t_0,T);\R^{d\times d})$ and $B(t,\cdot)\equiv B(t)\in L^\infty((t_0,T);\R^{d\times k})$, then $\cN_{\tau_1}(u)$ is congruent to the controllability Gramian $W_c$ of the LTV system $\dot{x}(t)=A(t)\,x(t)+B(t)\,u(t)$, and $\cN_{\tau_2}(u)=W_c$.  

(2) For any $y\in\R^d$,
\begin{equation}\label{eq:on the invertibility of N_1 and N_2}
y^\top \cN_{\tau_i}(u)\,y=\int_{t_0}^T |y^\top D\Phi_{t,\tau_i}(x_u(t))\,B(t,x_u(t))|^2\,dt, 
\end{equation}
hence $\cN_{\tau_i}(u)$ is symmetric and positive semi-definite (PSD).

(3) The Gramian map\footnote{ We use $\cS_d^+(\R)$ to denote the set of symmetric positive semi-definite matrices of order $d$ with real coefficients.} $\cN_{\tau_i}:L^\infty((t_0, T); \R^k)\to\cS_d^+(\R)$ defined by~\eqref{eq:nonautonomous N_i} is in general different from the notion of Gramian used in the homotopy continuation method \cite{chitour2006continuation,zhengping2024regularized,sussmann1992new} where the Gramian is defined as
\begin{equation}\label{eq:HCM gramian}
     M(u):=\int_{t_0}^{T}\,R_u(T,t)\,B(t,x_u(t))\,B(t,x_u(t))^\top\,R_u(T,t)^\top\,dt,
\qquad\, u\in L^\infty((t_0, T); \R^k)
\end{equation}
where $R_u(t,s)$ is defined in Lemma~\ref{lem:state-transition matrix of the nonlinear control B=0}.
\end{remark}

The proof of the following proposition uses standard arguments of minimization in Hilbert spaces. 
\begin{proposition}\label{pro:min-norm-solutions}
Let $x^0\in\R^d$, $u\in L^\infty((t_0, T); \R^k)$, and $x_u\in C^0([t_0,T]; \R^d)$ be the solution to~\eqref{eq:state-dependent-input}. Fix $i\in\{1,2\}$ and define the bounded linear operator $L_{u,\tau_i}:L^2((t_0, T);\R^k)\to\R^d$ by
\begin{equation}\label{eq:map L_u tau_i}
    L_{u,\tau_i}v = \int_{t_0}^TD\Phi_{t,\tau_i}(x_{u}(t))\,B(t,x_u(t))\,v(t)\,dt,\qquad\tau_1:=t_0,\,\tau_2:=T.
\end{equation}
For $y_i\in\Ran(L_{u,\tau_i})$, the problem
\begin{equation}\label{eq:minimizations problems}
\min\{\|v\|_2: L_{u,\tau_i}v=y_i\}
\end{equation}
has the solution $v_i(u):=v_i(u,x^0,y_i)$ given by
\begin{equation}\label{eq:relative minimal control}
v_i(u)(t) = \left[L_{u,\tau_i}^\ast\,\cN_{\tau_i}(u)^\dagger\,y_i\right](t) = B(t,x_u(t))^\top D\Phi_{t,\tau_i}(x_{u}(t))^\top\,\cN_{\tau_i}(u)^\dagger\, y_i,\qquad 
t\in[t_0, T].
\end{equation}
Here $L_{u,\tau_i}^\ast$ is the adjoint of $L_{u,\tau_i}$, and $\cN_{\tau_i}(u)^\dagger$ denotes the Moore--Penrose pseudoinverse of $\cN_{\tau_i}(u)$.
\end{proposition}

Proposition~\ref{pro:min-norm-solutions} ensures that $v_i(u)$ is $L^2$–minimal within the affine subspace defined by $L_{u,\tau_i}$.  
In the nonlinear setting, these subspaces depend on $u$, so minimality is relative rather than global. 

Using~\eqref{eq:EGectral norm of D_s_t_phi} and the singular value decomposition of symmetric PSD matrices, one finds that
\begin{equation}\label{eq:bound on supremum of v_i}
    \|v_i(u)\|_\infty\le\frac{\|B\|_\infty e^{\Lambda_1\Delta t_0}|y_i|}{\lambda_{\min}^\ast(\cN_{\tau_i}(u))}\qquad\forall u\in L^\infty((t_0, T); \R^k),\ i\in\{1,2\},
\end{equation}
where $\lambda_{\min}^\ast(\cN_{\tau_i}(u))>0$ is the smallest nonzero eigenvalue of $\cN_{\tau_i}(u)$.

\paragraph{Feasible coercivity class}
Let $\lambda_{\min}\big(\cN_{\tau_i}(u)\big)$ denote the smallest eigenvalue of $\cN_{\tau_i}(u)$, and let $C_i>0$ be a constant depending only on system data (e.g.\ $\|B\|_\infty$, $\Lambda_1$, $\Lambda_2$, $\Delta t_0$, and possibly $x^0$).  Guided by~\eqref{eq:bound on supremum of v_i}, we encode the positivity of $\cN_{\tau_i}(u)$ through the \emph{feasible coercivity class}
\begin{equation}\label{eq:feasible coercivity class}
    \cF(C_i):=\Big\{u\in L^\infty((t_0, T);\R^k):\ \lambda_{\min}\big(\cN_{\tau_i}(u)\big)\ge C_i^{-1}\Big\},\quad i\in\{1,\,2\}.
\end{equation}

The following observations on $\cF(C_i)$ are immediate.  

\begin{remark}\label{rmk:UCi-properties}
\begin{enumerate}
\item \emph{Possibility of emptiness.} For a given $C_i>0$, the feasible coercivity class $\cF(C_i)$ may be empty. Non-emptiness requires the existence of at least one control $u$ with $\lambda_{\min}(\cN_{\tau_i}(u))\ge C_i^{-1}$.
\item \emph{Monotonicity in $C_i$.} If $C_i\le C_i'$, then $\cF(C_i)\subseteq \cF(C_i')$.
\item \emph{Topological property.} By Assumptions~\ref{ass:general assumption nonautonomous}, the map $u\mapsto \cN_{\tau_i}(u)$ is continuous, hence $u\mapsto \lambda_{\min}(\cN_{\tau_i}(u))$ is continuous. Therefore $\cF(C_i)$ is closed in $L^\infty((t_0,T);\R^k)$. We clarify this point later in Lemma~\ref{lem:Lipschitz constant of N_i} and Corollary~\ref{cor:N_i continuity}.
\item \emph{Non-emptiness criterion.} If there exists a control reference $u^\flat\in L^\infty((t_0,T);\R^k)$ with $\lambda_{\min}(\cN_{\tau_i}(u^\flat))\ge C_i^{-1}$, then $\cF(C_i)\neq\emptyset$. By the continuity of $u\mapsto\lambda_{\min}(\cN_{\tau_i}(u))$, a neighborhood of $u^\flat$ also lies in $\cF(C_i)$. A sufficient condition ensuring $\cF(C_i)=L^\infty((t_0,T);\R^k)$ for every $x^0\in\R^d$ is provided in Proposition~\ref{pro:uniformly invertible B}.
\end{enumerate}
\end{remark}

For any $u\in\cF(C_i)$, estimate~\eqref{eq:bound on supremum of v_i} yields the uniform bound
\begin{equation}\label{eq:global-uniform-bound}
\|v_i(u)\|_\infty\;\le\zeta(y_i)= C_i\,\|B\|_\infty\,e^{\Lambda_1\Delta t_0}\,|y_i|,\qquad i\in\{1,2\}.
\end{equation}
Let $v_i(u,x^0,y_i)$ is defined as in Proposition~\ref{pro:min-norm-solutions}, we define the synthesis maps
\begin{equation}\label{eq:synthesis maps}
   \cS_i:L^\infty((t_0, T); \R^k)\to L^\infty((t_0, T); \R^k),\;u\mapsto v_i(u,x^0,y_i),\qquad y_i\in\R^d ,\quad i\in\{1,\,2\}.
\end{equation}
Motivated by~\eqref{eq:global-uniform-bound}, we introduce the feasibility coercivity ball
\begin{equation}\label{eq:feasibility ball}
    \cF(y_i):=\cF(C_i)\cap\cB_i\quad\text{where}\quad\cB_i:=\{u\in L^\infty((t_0, T); \R^k):\|u\|_\infty\le\,\zeta_i:=\zeta(y_i)\}.
\end{equation}

Under the non-emptiness assumption on $\cF(C_i)$, and if $\cS_i$ acts as a self-map on $\cF(C_i)$, namely
\begin{equation}\label{eq:self-mapping}
    \cS_i(\cF(C_i))\subset\cF(C_i),
\end{equation}
we can prove the existence and uniqueness of fixed points for $\cS_i$ in $\cF(y_i)$, whenever an associated Volterra rank-one comparison estimate is contractive.

The proof of the following result is given in Section~\ref{ss:abstract-fixed-point}. 

\begin{theorem}\label{thm:abstract-fixed-point}
Fix $i\in\{1,2\}$. Assume that $\cF(C_i)\neq\emptyset$ and that~\eqref{eq:self-mapping} is satisfied. Then there exists a positive sequence $(\varrho_m)_{m\ge1}$ with $\varrho_m:=K_i^m\gamma_m(p_i)$ such that
\begin{equation}\label{eq:supnorm of F_i^m(u)-F_i^m(v)}
    \|\cS_i^{\,m}(u)-\cS_i^{\,m}(v)\|_\infty
    \le \varrho_m\|u-v\|_\infty,
    \qquad \forall\,u,v\in\cF(y_i),\ \forall\,m\ge1,
\end{equation}
where $K_i>0$ depends only on system data and $C_i$, and for some $p_i\in(0,1)$ and $\gamma_m(p_i)\to 0$ as $m\to\infty$. Furthermore, if there exists $m_i\ge1$ such that 
\begin{equation}\label{eq:eventual-contraction-condition}
    \varrho_{m_i}<1,
\end{equation}
then $\cS_i$ admits a unique fixed point $u_i\in\cF(y_i)$, and the Picard iteration
$u^{(n+1)}=\cS_i(u^{(n)})$ converges to $u_i$ for any $u^{(0)}\in\cF(y_i)$.
\end{theorem}
\begin{remark}
   If $\cF(C_i)\neq\emptyset$ and~\eqref{eq:self-mapping} is satisfied, then $\cF(y_i)\ne\emptyset$. In fact, under these assumptions, $\cS_i(u)\in\cF(y_i)$ for all $u\in\cF(C_i)$. Furthermore, if $B(t,x)\equiv B(t)$ is not state-dependent, the decay estimate~\eqref{eq:supnorm of F_i^m(u)-F_i^m(v)} holds for all $(u, v)\in\cF(C_i)^2$ since the Gramian map $\cN_{\tau_i}$ is Lipschitz continuous on the whole $L^\infty((t_0, T); \R^k)$.
\end{remark}

Theorem~\ref{thm:abstract-fixed-point} guarantees the existence and uniqueness of a fixed point of the synthesis map $\cS_i$ on $\cF(y_i)$.
To connect this fixed point with $L^2$–norm minimality, we invoke a Lagrange multiplier. We will need preparatory results stated in Lemma~\ref{lem:preparatory result} below. To this end, we introduce for any fixed $x^0\in\R^d$ the endpoint map 
\begin{equation}\label{eq:endpoint map}
    \cE:=\cE_{x^0,T}:u\in\,L^\infty((t_0, T); \R^k)\longmapsto\R^d,\quad\cE(u) = x_u(T)
\end{equation}
and for $i\in\{1,2\}$ and any fixed $y_i\in\R^d$, the feasible map
\begin{equation}\label{eq:feasible map}
    G_{\tau_i}:L^\infty((t_0,T); \R^k)\to\R^d,\quad G_{\tau_i}(u) = L_{u,\tau_i}u-y_i
\end{equation}
where $x_u(\cdot)$ is the solution of~\eqref{eq:state-dependent-input} and $L_{u,\tau_i}$ is defined in~\eqref{eq:map L_u tau_i}. The proof of the following is given in Section~\ref{ss:preparatory result}.

\begin{lemma}\label{lem:preparatory result}
    One has $\cE,\,G_{\tau_i}\in C^1(\cY; \R^d)$. Furthermore, let $D\cE(u),\,DG_{\tau_i}(u):L^2((t_0,T);\R^k)\to\R^d$ are, respectively the Fréchet derivative of $\cE$ and $G_{\tau_i}$ at any fixed $u\in\cY$. Then, one has the following identity
  \begin{equation}\label{eq:identities}
     DG_{\tau_i}(u)= D\Phi_{T,\tau_i}\!\big(x_u(T)\big)\,D\cE(u).
  \end{equation}
  
 Moreover, let $\cN_{\tau_i}(u)$ and $M(u):=D\cE(u)\,D\cE(u)^\ast$ are respectively defined in~\eqref{eq:nonautonomous N_i} and~\eqref{eq:HCM gramian}.
  \begin{enumerate}
      \item If $\cN_{\tau_i}(u)$ is invertible and $\idty+K_{u,\tau_i}L_{u,\tau_i}^\ast\cN_{\tau_i}(u)^{-1}\in\R^{d\times d}$
      is invertible,
      then $DG_{\tau_i}(u)L_{u,\tau_i}^\ast\in\R^{d\times d}$ is invertible and $DG_{\tau_i}(u)$ is right-invertible;
      \item If $M(u)$ is invertible and $\idty-K_{u,\tau_i}DG_{\tau_i}(u)^\ast M(u)^{-1}\in\R^{d\times d}$ is invertible,
      then $L_{u,\tau_i}DG_{\tau_i}(u)^\ast\in\R^{d\times d}$ is invertible and $L_{u,\tau_i}$ is right-invertible.
  \end{enumerate}
 Here, we let $K_{u,\tau_i}:=DG_{\tau_i}(u)-L_{u,\tau_i}$.
\end{lemma}

The following results connect the fixed-point of Theorem~\ref{thm:abstract-fixed-point} with $L^2$-nom minimality in the feasible set; its proof is presented in Section~\ref{ss:minimizer-is-fp}.

\begin{theorem}\label{thm:minimizer-is-fp}
Fix $i\in\{1,2\}$, assume $\cF(C_i)\neq\emptyset$, and  $\idty+K_{u,\tau_i}L_{u,\tau_i}^\ast\cN_{\tau_i}(u)^{-1}\in\R^{d\times d}$ is invertible for all $u\in\cF(C_i)$. Let $y_i\in\R^d$ and the feasible set
\[
\mathfrak F_i:=\{u\in\cF(y_i):\ G_{\tau_i}(u)=0\}.
\]
Then, the following are equivalent for any local minimizer $\bar u$ of $\tfrac12\|u\|_{L^2}^2$ over $\mathfrak F_i$:

1. $\bar u$ is a fixed point of $\cS_i$ on $\cF(y_i)$, i.e., $\bar u=L_{\bar u,\tau_i}^\ast\,\cN_{\tau_i}(\bar u)^{-1}\,y_i$;

2. The following orthogonality condition holds 
\begin{equation}\label{eq:orthogonality condition}
  \langle[L_{\bar u,\tau_i}DG_{\tau_i}(\bar u)^\ast]^{-1}y_i,DG_{\tau_i}(\bar u)\,h\rangle_{\R^d}=0\qquad\forall\,h\in\ker L_{\bar u,\tau_i}.
\end{equation}
\end{theorem}

The following shows the existence of a local minimizer of $\tfrac12\|u\|_{L^2}^2$ over $\mathfrak F_i$. The proof is given in Section~\ref{ss:exist-on-ball-global}.

\begin{proposition}\label{prop:exist-on-ball-global}
 Fix $i\in\{1,2\}$ and assume $\cF(C_i)\neq\emptyset$. Let $y_i\in\R^d$, then the feasible set
 \begin{equation}
   \mathfrak F_i:=\{u\in\cF(y_i):\ G_{\tau_i}(u)=0\}
\end{equation}
 is weakly sequential closed in $L^2((t_0, T); \R^k)$. Consequently, $\tfrac12\|u\|_{L^2}^2$ attains a minimum on $\mathfrak F_i$.
\end{proposition}

From Theorems~\ref{thm:abstract-fixed-point}, and~\ref{thm:minimizer-is-fp}, and Proposition~\ref{prop:exist-on-ball-global}, we deduce the synthesis and optimal control results:

\begin{theorem}\label{thm:main controllability result nonlinear}
Let $(x^0,\,x^1)\in\R^d\times\R^d$, $i\in\{1,2\}$, and set $y_i:=\Phi_{T,\tau_i}(x^1)-\Phi_{t_0,\tau_i}(x^0)$.
Assume $\cF(C_i)\neq\emptyset$, \eqref{eq:self-mapping} and~\eqref{eq:eventual-contraction-condition} are satisfied, and let $u_i\in\cF(y_i)$ be the fixed point of $\cS_i$.
Then $u_i$ steers~\eqref{eq:state-dependent-input} from $x^0$ to $x^1$ on $[t_0,T]$ and
\begin{equation}\label{eq:control with N_i}
u_i(t)=L_{u_i,\tau_i}^\ast\,\cN_{\tau_i}(u_i)^{-1}y_i=B(t, x_{u_i}(t))^\top D\Phi_{t,\tau_i}(x_{u_i}(t))^\top\,\cN_{\tau_i}(u_i)^{-1}y_i,\quad t\in[t_0,T],\quad \tau_1=t_0,\, \tau_2=T,
\end{equation}
with $\|u_i\|_{L^2}^2 = y_i^\top\,\cN_{\tau_i}(u_i)^{-1}y_i$.

If the hypotheses of Theorem~\ref{thm:minimizer-is-fp} hold for every local minimizer over $\mathfrak F_i:=\{u\in\cF(y_i):L_{u,\tau_i}u=y_i\}$, and if~\eqref{eq:orthogonality condition} holds whenever needed, then $\tfrac12\|u\|_{L^2}^2$ attains its minimum over $\mathfrak F_i$, and the unique minimizer is $u_i$.
Equivalently,
\[
\|u_i\|_{L^2}^2\le \|w\|_{L^2}^2\quad\forall\,w\in\mathfrak F_i,\qquad \text{with equality iff }w=u_i.
\]
\end{theorem}

\begin{proof}
Since $u_i=\cS_i(u_i)=L_{u_i,\tau_i}^\ast\,\cN_{\tau_i}(u_i)^{-1}y_i$, we obtain \eqref{eq:control with N_i}, feasibility $L_{u_i,\tau_i}u_i=y_i$, and $\|u_i\|_{L^2}^2=y_i^\top\cN_{\tau_i}(u_i)^{-1}y_i$.
By Proposition~\ref{prop:exist-on-ball-global}, $\tfrac12\|u\|_{L^2}^2$ attains a minimum on $\mathfrak F_i$. Let $\bar u$ be a minimizer. Under the stated hypotheses, Theorem~\ref{thm:minimizer-is-fp} implies that $\bar u$ is a fixed point of $\cS_i$. By uniqueness, $\bar u=u_i$.
\end{proof}

\begin{remark}[A posteriori check of the minimizer]
When $\cS_i$ admits a unique fixed
point $u_i\in\cF(y_i)$, it suffices to verify whether the orthogonality condition~\eqref{eq:orthogonality condition} holds at $u_i$, viz.
\begin{equation}\label{eq:a posteriori RAC}
 \langle[L_{u_i,\tau_i}DG_{\tau_i}(u_i)^\ast]^{-1}y_i,DG_{\tau_i}(u_i)\,h\rangle_{\R^d}=0\qquad\forall\,h\in\ker L_{u_i,\tau_i}.    
\end{equation}
In this case, Theorem~\ref{thm:minimizer-is-fp} implies that $u_i$ is a local minimizer of $\frac12\|u\|_{L^2}^2$ over $\mathfrak F_i$, and that if the
global minimizer $\bar u\in\mathfrak F_i$ satisfies~\eqref{eq:orthogonality condition}, then it coincides with $u_i$; hence $u_i$ is the
global minimum–energy control on $\mathfrak F_i$. 
\end{remark}

\begin{remark}
	If {$B(t,x)\equiv B(t)$} and there exists $C_i>0$ such that $\cF(C_i)=L^\infty((t_0,T);\R^k)$, then the conclusions of Theorem~\ref{thm:main controllability result nonlinear} hold on  {any bounded set $\cB\subset L^\infty$ containing $\cF(y_i)$}. In particular, if~\eqref{eq:a posteriori RAC} holds {on $\cB$}, the fixed point $u_i$ is the unique global minimizer of $\tfrac12\|u\|_{L^2}^2$ among all {$u\in\cB$} satisfying $L_{u,\tau_i}u=y_i$.
\end{remark}

\begin{remark}[Relation to reachability for fixed $T>t_0$]\label{rmk:relation to reachability} 
Theorem~\ref{thm:main controllability result nonlinear} is fundamentally a \emph{synthesis result}: given the fixed point $u_i\in\cF(y_i)$ of $\cS_i$, it provides an explicit representation of a control that steers $x^0$ to $x^1$. 
It does not assert that such a fixed point exists for every pair $(x^0,x^1)$, nor that the reachable set coincides with $\R^d$. 
Reachability depends on the structure of the feasible coercivity class $\cF(C_i)$ and on whether the self-mapping property~\eqref{eq:self-mapping} is satisfied. 
Three distinct situations may occur:

\emph{(a) Complete controllability.}  
If there exists $C_i>0$, independent of the initial state $x^0\in\R^d$, such that $\cF(C_i)=L^\infty((t_0,T);\R^k)$, then~\eqref{eq:self-mapping} holds automatically and $\cR_\Sigma(T,x^0)=\R^d$ for all $x^0\in\R^d$. 
In particular, system~\eqref{eq:state-dependent-input} is completely controllable on $[t_0, T]$. The converse is not true: complete controllability does not imply the existence of $C_i>0$ such that $\cF(C_i)=L^\infty((t_0, T);\R^k)$. 

\emph{(b) Controllability from a fixed $x^0$.}  
If for some $C_i=C_i(x^0)>0$ depending on the initial condition $x^0$, one has $\cF(C_i)=L^\infty((t_0,T);\R^k)$, then~\eqref{eq:self-mapping} again holds, and $\cR_\Sigma(T,x^0)=\R^d$. 
Thus, the system is controllable from that specific $x^0$, although the property may not extend to all other initial states.

\emph{(c) Restricted synthesis.}  
If for some $C_i=C_i(x^0)>0$ one only has $\emptyset\neq\cF(C_i)\subsetneq L^\infty((t_0,T);\R^k)$, then Theorem~\ref{thm:main controllability result nonlinear} ensures that any fixed point $u_i\in\cF(y_i)$ steers $x^0$ to $x^1$, provided~\eqref{eq:self-mapping} is verified. 
Thus, the synthesis produces a family of feasibly controls, but the reachable set $\cR_\Sigma(T,x^0)$ may be a proper subset of $\R^d$; its precise description depends on the estimates involving $|\Phi_{T,t_0}(x^1)-x^0|$ or $|x^1-\Phi_{t_0,T}(x^0)|$ that guarantee~\eqref{eq:self-mapping}.
\end{remark}

Directly checking invertibility of $\cN_{\tau_i}(u)$ is generally a difficult task since it depends on the full trajectory $x_u(\cdot)$. 
This motivates the search for structural assumptions under which, for some $C_i>0$, $\cF(C_i)\neq\emptyset$ or coincides with $L^\infty((t_0, T);\R^k)$. The first such condition is presented below. The proof is presented in Section~\ref{ss:uniformly invertible B}.

\begin{proposition}\label{pro:uniformly invertible B}
Fix $i\in\{1,2\}$. Suppose $d=k$ and there exists a nonzero
$b\in L^1((t_0,T);\R_+)$ such that
\begin{equation}\label{eq:not pointwise invertibility of B(t)}
    |B(t,x)^\top y|\ \ge\ b(t)\,|y|\qquad\text{for a.e. } t\in(t_0,T),\,\forall (x, y)\in\R^d\times\R^d.
\end{equation}
Then $\cN_{\tau_i}(u):=\cN_{\tau_i}(u,x^0)$ is uniformly coercive, viz.
\begin{equation}\label{eq:coercivity bound full actuated case}
    y^\top\,\cN_{\tau_i}(u)\,y\ \ge\ \frac{e^{-2\Lambda_1\Delta t_0}\,\|b\|_1^2}{\Delta t_0}\,|y|^2\qquad \forall (y,\,x^0,\,u)\in \R^d\times\R^d\times L^\infty((t_0,T);\R^k).
\end{equation}
Hence, one may take $C_i=\Delta t_0\,e^{2\Lambda_1\Delta t_0}\,\|b\|_1^{-2}$, so that $\cF(C_i)=L^\infty((t_0,T);\R^k)$. In particular, \eqref{eq:self-mapping} is satisfied.
\end{proposition}

\medskip
Proposition~\ref{pro:uniformly invertible B} addresses the fully actuated case $k=d$ under a structural assumption on $B(t,x)$ ensuring uniform coercivity. 
If this assumption fails, or more generally in the underactuated regime $k<d$, such global bounds are no longer guaranteed. The following analysis provides sufficient conditions for reachability and synthesis in this configuration. The proof is provided in Section~\ref{ss:relative coercivity from a control reference}.

\begin{theorem}\label{thm:relative coercivity from a control reference}
    Let $x^0\in\R^d$ and fix $i\in\{1,2\}$. Assume the existence of a reference control $u_i^\flat:=u_i^\flat(x^0)\in L^\infty((t_0, T); \R^k)$ such that $\cN_{\tau_i}(u_i^\flat)$ is invertible. Let $C_i:=(1+\theta)/\lambda_{\min}\big(\cN_{\tau_i}(u_i^\flat)\big)$ where $\theta\in(0,1]$.
    Furthermore, assume that the following estimate holds
    \begin{equation}\label{eq:bound on the sup-norm of the reference control}
        \|u_i^\flat\|_\infty< \frac{\theta\,\lambda_{\min}(\cN_{\tau_i}(u_i^\flat))}{(1+\theta)\,L_{\cN_{\tau_i}}}.
    \end{equation}
    Then $\cF(C_i)\neq\emptyset$, and~\eqref{eq:self-mapping} holds if $y_i\in\R^d$ satisfies
    \begin{equation}\label{eq:sufficient conditions for invertibility general}
        |y_i|\le  \vartheta(u_i^\flat)\ :=\ \frac{\lambda_{\min}(\cN_{\tau_i}(u_i^\flat))}{(1+\theta)\,\|B\|_\infty\,e^{\Lambda_1\Delta t_0}}
\ \left(\frac{\theta\,\lambda_{\min}(\cN_{\tau_i}(u_i^\flat))}{(1+\theta)\,L_{\cN_{\tau_i}}}-\|u_i^\flat\|_\infty\right). 
    \end{equation}
    Here $L_{\cN_{\tau_i}}>0$ is the Lipschitz constant of the map $u\in\{u\in\,L^\infty((t_0, T); \R^k):\|u\|_\infty\le\,\zeta_i \}\mapsto\cN_{\tau_i}(u)$.
\end{theorem}

In contrast to Proposition~\ref{pro:uniformly invertible B}, Theorem~\ref{thm:relative coercivity from a control reference} applies under Assumption~\ref{ass:general assumption nonautonomous} alone, requiring no additional structure on $N_t$ or $B(t,x)$. 
Before discussing its implications, we first establish the Lipschitz constant $L_{\cN_{\tau_i}}>0$ for $\cN_{\tau_i}(\cdot)$ under Assumptions~\ref{ass:general assumption nonautonomous} and~\ref{ass:gneneral assumption on B}. The proof of the following result is given in Section~\ref{ss:proof of Lipschitz constant of N_i}.

\begin{lemma}\label{lem:Lipschitz constant of N_i}
    Fix $i\in\{1,2\}$. Then $\cN_{\tau_i}$ defined by~\eqref{eq:nonautonomous N_i} is continuous from $L^\infty((t_0, T); \R^k)$ to $\cS_d^+(\R)$ and
    globally Lipschitz from $\cB_i:=\{u\in\,L^\infty((t_0, T); \R^k):\|u\|_\infty\le\,\zeta_i \}$ to $\cS_d^+(\R)$. For all $(u,\,v)\in\cB_i^2$,
        \begin{align}\label{eq:Lipschitz constant of N_1}
\frac{\|\cN_{\tau_1}(u)-\cN_{\tau_1}(v)\|}{\|u-v\|_\infty}
&\le\,L_{\cN_{\tau_1}}:=\frac{L_B\|B\|_\infty^2 e^{\Lambda_1\Delta t_0}}{\Lambda_1}
\Bigg[\frac{e^{(2\Lambda_1+L_B\zeta_1)\Delta t_0}-e^{\Lambda_1\Delta t_0}}{\Lambda_1+L_B\zeta_1}
      +\frac{e^{(2\Lambda_1+L_B\zeta_1)\Delta t_0}-e^{-\Lambda_1\Delta t_0}}{L_B\zeta_1+3\Lambda_1}\Bigg]\nonumber \\
&\quad\qquad + \frac{2\|B\|_\infty^3\Lambda_2 e^{\Lambda_1\Delta t_0}}{\Lambda_1(\Lambda_1+L_B \zeta_1)}
\Bigg[
\frac{e^{(3\Lambda_1+L_B\zeta_1)\Delta t_0}-e^{-\Lambda_1\Delta t_0}}{4\Lambda_1+L_B\zeta_1}
-\frac{e^{(2\Lambda_1+L_B\zeta_1)\Delta t_0}-e^{-\Lambda_1\Delta t_0}}{3\Lambda_1+L_B\zeta_1} \nonumber\\
&\qquad\qquad\qquad\qquad\qquad\qquad
+\frac{e^{\Lambda_1\Delta t_0}-e^{-\Lambda_1\Delta t_0}}{2\Lambda_1}
-\frac{e^{2\Lambda_1\Delta t_0}-e^{-\Lambda_1\Delta t_0}}{3\Lambda_1}
\Bigg].
\end{align}
    \begin{align}\label{eq:Lipschitz constant of N_2}
\frac{\|\cN_{\tau_2}(u)-\cN_{\tau_2}(v)\|}{\|u-v\|_\infty}
&\le\,L_{\cN_{\tau_2}}:=\frac{L_B\|B\|_\infty^2 e^{\Lambda_1\Delta t_0}}{\Lambda_1}
\Bigg[\frac{e^{\Lambda_1\Delta t_0}-e^{L_B\zeta_2\Delta t_0}}{\Lambda_1-L_B\zeta_2}
      -\frac{e^{L_B\zeta_2\Delta t_0}-e^{-\Lambda_1\Delta t_0}}{L_B\zeta_2+\Lambda_1}\Bigg]\nonumber \\
&\quad\qquad\qquad + \frac{2\|B\|_\infty^3\Lambda_2 e^{\Lambda_1\Delta t_0}}{\Lambda_1(\Lambda_1+L_B \zeta_2)}
\Bigg[
\frac{e^{2\Lambda_1\Delta t_0}-e^{L_B\zeta_2\Delta t_0}}{2\Lambda_1-L_B\zeta_2}
+\frac{e^{\Lambda_1\Delta t_0}-e^{L_B\zeta_2\Delta t_0}}{\Lambda_1-L_B\zeta_2} \nonumber\\
&\qquad\qquad\qquad\qquad\qquad\qquad\qquad
+\frac{e^{\Lambda_1\Delta t_0}-e^{-\Lambda_1\Delta t_0}}{2\Lambda_1}
-\frac{e^{2\Lambda_1\Delta t_0}-e^{-\Lambda_1\Delta t_0}}{3\Lambda_1}
\Bigg].
\end{align}
\end{lemma}
\begin{remark}\label{rmk:specific case of Lipschitz constants}
    When $B(t,x)\equiv\,B(t)$, \eqref{eq:Lipschitz constant of N_1} and~\eqref{eq:Lipschitz constant of N_2} holds for $u,\,v\in L^\infty((t_0, T); \R^k)$, and
    \[
    L_{\cN_{\tau_1}}:=\frac{\Lambda_2\,\|B\|_\infty^3}{6}\left(3e^{\Lambda_1\Delta t_0}+1\right)\left(\frac{e^{\Lambda_1\Delta t_0}-1}{\Lambda_1}\right)^3,\quad\,L_{\cN_{\tau_2}}:=\frac{\Lambda_2\,\|B\|_\infty^3}{3}\left(\frac{e^{\Lambda_1\Delta t_0}-1}{\Lambda_1}\right)^3.
    \]
\end{remark}
The following result is then an immediate consequence.
\begin{corollary}\label{cor:N_i continuity}
By Lemma~\ref{lem:Lipschitz constant of N_i} and spectral stability
\[
|\lambda_{\min}(C)-\lambda_{\min}(D)|\le\|C-D\|\qquad\forall (C,\,D)\in\cS_d(\R)^2,
\]
the map $u\ \mapsto\ \lambda_{\min}\!\big(\cN_{\tau_i}(u)\big)$
is continuous on $L^\infty((t_0, T); \R^k)$ and Lipschitz on $\{u\in\,L^\infty((t_0, T); \R^k):\|u\|_\infty\le\,\zeta_i \}$. In particular, the feasible coercivity class $\cF(C_i)$ defined in~\eqref{eq:feasible coercivity class}
is closed in $L^\infty((t_0, T); \R^k)$. 
\end{corollary}

\begin{remark}\label{rmk:on the Lipschitz constant of cN_i}
The bounds in~\eqref{eq:sufficient conditions for invertibility general}, \eqref{eq:Lipschitz constant of N_1}, and~\eqref{eq:Lipschitz constant of N_2}, while ``sharp'' in a general sense (under Assumption~\ref{ass:general assumption nonautonomous} alone), are conservative for specific vector fields $N_t$. 
In particular, the exponential factor $e^{\Lambda_1\Delta t_0}$ reflects the global $\Lambda_1$–Lipschitz bound on $N_t$ (cf.\ Lemma~\ref{lem:spectral norm of D_s_t_phi and D^2_s_t_phi}), whereas for certain dynamics the differential $\|D\Phi_{t,s}(\cdot)\|$ can in fact decay exponentially. We provide illustrative examples later in Section~\ref{ss:refined estimates for HNN}. 
\end{remark}

\subsubsection{A practical control reference: the zero control}
 In practice, to check the assumptions of Theorem~\ref{thm:relative coercivity from a control reference}, a natural first choice is $u_i^\flat\equiv0=:u_0$, for which the corresponding trajectory reads
\[
x_{u_0}(t)=\Phi_{t_0,t}(x^0),\qquad t\in[t_0,T].
\]
For $i\in\{1,2\}$ and $\tau_1=t_0,\ \tau_2=T$, this yields the symmetric positive semi-definite matrices
\begin{equation}\label{eq:SSD W_i}
W_i(T):=\cN_{\tau_i}(0, x^0)
=\int_{t_0}^{T}\! D\Phi_{t,\tau_i}\!\big(\Phi_{t_0,t}(x^0)\big)\,B(t,\Phi_{t_0,t}(x^0))B(t,\Phi_{t_0,t}(x^0))^\top\,D\Phi_{t,\tau_i}\!\big(\Phi_{t_0,t}(x^0)\big)^\top\,dt.
\end{equation}
It is immediate from \eqref{eq:algebraic identities 2} and \eqref{eq:equation derivative of phi_t_0,t} that $W_2(T)\in\R^{d\times d}$ solves the time-varying Lyapunov differential equation
\begin{equation}\label{eq:evolution of W_2}
\dot W(t)=B(t,\Phi_{t_0,t}(x^0))B(t,\Phi_{t_0,t}(x^0))^\top+DN_t\!\big(\Phi_{t_0,t}(x^0)\big)\,W(t)+W(t)\,DN_t\!\big(\Phi_{t_0,t}(x^0)\big)^\top,\quad
W(t_0)=0.
\end{equation}
Moreover, by Remark~\ref{rmk:matrix  S_{t_0,x^0}(s,t)} one has,
\begin{equation}\label{eq:relation between W_1 and W_2}
W_2(t)=S_{t_0,x^0}(t,t_0)\,W_1(t)\,S_{t_0,x^0}(t,t_0)^\top,\qquad t\in[t_0,T]
\end{equation}
so $W_1(t)$ is invertible iff $W_2(t)$ is invertible.

The following result is then an immediate consequence of Theorem~\ref{thm:relative coercivity from a control reference}.
\begin{corollary}\label{cor:relative coercivity from a zero control reference}
    Let $x^0\in\R^d$, $T>t_0$, and fix $i\in\{1,2\}$. Assume that $W_2(T)$ is invertible, and let $C_i:=(1+\theta)/\lambda_{\min}(W_i(T))$ where $\theta\in(0,1]$. Then $\cF(C_i)\neq\emptyset$, and~\eqref{eq:self-mapping} holds if $y_i\in\R^d$ satisfies
    \begin{equation}\label{eq:sufficient conditions for invertibility zero control}
        |y_i|\le\frac{\theta\,\lambda_{\min}^2(W_i(T))}{(1+\theta)^2\,L_{\cN_{\tau_i}}\,\|B\|_\infty\,e^{\Lambda_1\Delta t_0}}.
    \end{equation}
    Here $L_{\cN_{\tau_i}}>0$ is the Lipschitz constant of the map $u\in\{u\in\,L^\infty((t_0, T); \R^k):\|u\|_\infty\le\,\zeta_i \}\mapsto\cN_{\tau_i}(u)$. Furthermore, for any $y_i\in\R^d$ satisfying~\eqref{eq:sufficient conditions for invertibility zero control},
    the following energy estimates hold
    \begin{equation}\label{eq:energy of u_i for the zero reference control}
         \int_{t_0}^{T}|u_i(t)|^2\, dt\le\frac{(1+\theta)}{\lambda_{\min}(W_i(T))}|y_i|^2
    \end{equation}
    where $u_i\in\cF(y_i)$ is the fixed point of the synthesis map $\cS_i:\cF(y_i)\to\cF(y_i)$ defined by~\eqref{eq:synthesis maps}.
\end{corollary}
\begin{proof}
    Inequality~\eqref{eq:sufficient conditions for invertibility zero control} follows from~\eqref{eq:sufficient conditions for invertibility general}. Next ~\eqref{eq:control with N_i} and~\eqref{eq:nonautonomous N_i} yields
    \begin{equation}
    \int_{0}^{T}|u_i(t)|^2\, dt=\int_{t_0}^Tu_i(t)^\top\,u_i(t)\,dt=y_i^\top\cN_{\tau_i}(u_i)^{-1}\,y_i\le\frac{|y_i|^2}{\lambda_{\min}(\cN_{\tau_i}(u_i))}\le\frac{(1+\theta)}{\lambda_{\min}(W_i(T))}|y_i|^2,
    \end{equation}
  since $\lambda_{\min}(\cN_{\tau_i}(u_i))\ge \lambda_{\min}(W_i(T))/(1+\theta)$.
\end{proof}
\begin{remark}
It follows from Remark~\eqref{rmk:matrix  S_{t_0,x^0}(s,t)} that $W_2(T)=\cN_{\tau_2}(0,x^0)$ is invertible iff the LTV system
\begin{equation}\label{eq:LTV system}
\dot y(t)=DN_t\!\big(\Phi_{t_0,t}(x^0)\big)\,y(t)+B(t,\Phi_{t_0,t}(x^0))\,u(t),\qquad
y(t_0)=y^0\in\R^d,\quad t\in[t_0,T],
\end{equation}
is controllable with inputs $u\in L^1((t_0,T);\R^k)$. This is the Kalman criterion for LTV systems via the controllability Gramian \cite[Theorem.~1.16]{coron2007control}. Under additional regularity, invertibility of $W_2(T)$ can also be characterized by algebraic rank conditions; see, e.g., \cite[Theorem~1]{chang1965algebraic}, \cite[Theorem~1.18]{coron2007control} and~\cite[Proposition~3.5.15]{sontag2013mathematical}.
\end{remark}

\subsubsection{Optimizing the admissible radius: balancing coercivity and control size}\label{ss:optimizing the reachable set estimates} The trivial choice $u_i^\flat=0$ might not be the optimal one in practice to apply Theorem~\ref{thm:relative coercivity from a control reference} if $W_2(T)$ is invertible. In fact, it follows from the admissible target displacement estimate~\eqref{eq:sufficient conditions for invertibility general} that a good reference $u_i^\flat\in L^\infty((t_0, T); \R^k)$ is one that maximizes $\vartheta(u)$, which jointly increases the positive coercivity contribution and controls the penalty induced by $\|u\|_\infty$. The aim of this section is to illustrate how we can choose a reference control that maximizes $\vartheta(u)$.

Let $m\in\N$, $m\ge 1$, and $\{\varphi_j(t):t\in[t_0, T]\}_{1\le j\le m}$ be a chosen (e.g., with respect to the properties of the vector field $N_t$ and the input matrix $B$) $m$-linearly independent familly of functions in $L^\infty((t_0, T);\R^k)$. Fix $i\in\{1,2\}$, $\theta\in(0,1]$ and introduce the set
\begin{equation}
    \cU_{\rm ad} := \left\{u\in \operatorname{Span}\{\varphi_j\}:g(u):=\frac{\theta\,\lambda_{\min}(\cN_{\tau_i}(u))}{(1+\theta)\,L_{\cN_{\tau_i}}}-\|u\|_\infty\ge 0\right\}.
\end{equation}
\begin{proposition}\label{pro:optimal admissible radius}
 Let $x^0\in\R^d$, $T>t_0$, and fix $i\in\{1,2\}$. Assume that $W_2(T)$ is invertible, and let
 \[
\vartheta(u)\ :=\ \frac{\lambda_{\min}(\cN_{\tau_i}(u))}{(1+\theta)\,\|B\|_\infty\,e^{\Lambda_1\Delta t_0}}
\left(\frac{\theta\,\lambda_{\min}(\cN_{\tau_i}(u))}{(1+\theta)\,L_{\cN_{\tau_i}}}-\|u\|_\infty\right),\qquad u\in L^\infty((t_0,T);\R^k).
\]
 Then, the following program 
 \begin{equation}\label{eq:optimal reference control}
     \max_{u\in\cU_{\rm ad}} \ \vartheta(u)
 \end{equation}
 admits at least one solution $u_i^\flat\in\cU_{\rm ad}$. Furthermore, if $\vartheta(u_i^\flat)>0$ then $\lambda_{\min}(\cN_{\tau_i}(u_i^\flat))>0$, $\cF(C_i)\neq\emptyset$, and~\eqref{eq:self-mapping} holds if $y_i\in\R^d$ satisfies $|y_i|\le \vartheta(u_i^\flat)$ where $C_i:=(1+\theta)/\lambda_{\min}\big(\cN_{\tau_i}(u_i^\flat)\big)$.
\end{proposition}

\begin{proof}
    Since $W_2(T)$ is invertible, $0\in\cU_{\rm ad}$ so that $\cU_{\rm ad}\neq\emptyset$, and $\cU_{\rm ad}$ is closed since $\cU_{\rm ad}=g^{-1}([0, \infty))$ and $u\mapsto g(u)$ is continuous. Furthermore, since $\lambda_{\min}(\cN_{\tau_i}(u))\le\|\cN_{\tau_i}(u)\|\le\|B\|_\infty^2(e^{2\Lambda_1\Delta t_0}-1)/(2\Lambda_1)$,
    one deduce that any $u\in\cU_{\rm ad}$ satisfies
    \[
    \|u\|_\infty\le\frac{\theta\,\|B\|_\infty^2}{(1+\theta)\,L_{\cN_{\tau_i}}}\frac{e^{2\Lambda_1\Delta t_0}-1}{2\Lambda_1}, 
    \]
    so that $\cU_{\rm ad}$ is bounded, and therefore compact. Since $u\mapsto\,\vartheta(u)$ is continuous, it follows that~\eqref{eq:optimal reference control} has at least one solution $u_i^\flat\in\cU_{\rm ad}$ by the Weierstrass extreme value theorem. The last part of the proposition follows from Theorem~\ref{thm:relative coercivity from a control reference}, and the following equivalence
\[
\vartheta(u_i^\flat)>0\quad\Longleftrightarrow\quad\left(\|u_i^\flat\|_\infty<\frac{\theta\,\lambda_{\min}(\cN_{\tau_i}(u_i^\flat))}{(1+\theta)\,L_{\cN_{\tau_i}}}\quad\text{and}\quad\lambda_{\min}(\cN_{\tau_i}(u_i^\flat))>0\right).
\] 
\end{proof}

\begin{remark}
In Proposition~\ref{pro:optimal admissible radius}, the assumption that $W_2(T)$ is invertible provides the simplest guarantee that $\cU_{\rm ad}\neq\emptyset$. More generally, one may only require the existence of $u_i^\ast\in L^\infty((t_0,T);\R^k)$ such that $\cN_{\tau_i}(u_i^\ast)$ is invertible and $g(u_i^\ast)\ge0$. Verifying the latter condition, however, may be nontrivial in concrete models.
\end{remark}

\subsubsection{Refined estimates for Hopfield-type dynamics}\label{ss:refined estimates for HNN} 
We now return to a key point: the estimates for $\|D\Phi_{t,s}(\cdot)\|$ (Lemma~\ref{lem:spectral norm of D_s_t_phi and D^2_s_t_phi}) and for $L_{\cN_{\tau_i}}$ (Lemma~\ref{lem:Lipschitz constant of N_i}) are somehow ``sharp'' under Assumption~\ref{ass:general assumption nonautonomous}, but may be overly conservative for concrete models since they use the global estimates~\eqref{eq:EGectral norm of DN_t} and~\eqref{eq:EGectral norm of D^2N_t}.  We illustrate this fact in this section. Consider the vector field from Hopfield-type recurrent neural networks~\cite{hopfield1984neurons}, widely used in theoretical neuroscience and large-scale brain modeling, viz.
\begin{equation}\label{eq:Hopfield VF}
    N(x) = -D\,x+W \sigma(x),\qquad x\in\R^d,
\end{equation}
where $D\in\R^{d\times d}$ is diagonal and positive, $W\in\R^{d\times d}$ encodes connectivity, and $\sigma$ is the neural activation,
\[
\sigma(x):=(\sigma_1(x_1),\dots,\sigma_d(x_d))^\top,\qquad x:=(x_1,\,\cdots,\,x_d)^\top\in\R^d,
\]
with $\sigma_i(\cdot)$ is of sigmoid-type, say, $\sigma_i(s)=\tanh(s)$. Although $W=W(t)$ may vary analytically in time, e.g.\ to capture neuron–astrocyte interactions~\cite{gong2024astrocytes}, it is commonly assumed to be constant, a simplification that we adopt. As a sigmoid function, $\sigma_i$ is infinitely differentiable and globally bounded on $\R$ with all its successive derivatives, often $\sigma_i'(\cdot)\ge 0$. It follows that $N$ satisfies all the hypotheses of Assumption~\ref{ass:general assumption nonautonomous}. Since $N$ is autonomous, one assumes that $t_0=0$. In particular, its generated flow $\phi_t:=\Phi_{0,t}$, $t\in\R$, is a one parametric subgroup of $\operatorname{Diff}(\R^d)$, the group of all diffeomorphisms on $\R^d$. In this framework, for $i\in\{1,2\}$, one has
\begin{equation}\label{eq:N_i for Hopfield vector field}
    \cN_{\tau_i}(u):=\cN_{\tau_i}(u, x^0) = \int_0^T\,D\phi_{\tau_i-t}(x_u(t))\,B(t)\,B(t)^\top\,D\phi_{\tau_i-t}(x_u(t))^\top\,dt,\qquad\tau_1:=0,\,\tau_2:=T,
\end{equation}
where the input matrix $B(t,x)\equiv B(t)\in\,L^\infty((t_0, T); \R^{d\times k})$, and $x_u$ is the corresponding solution to~\eqref{eq:state-dependent-input} with the vector field in~\eqref{eq:Hopfield VF}. For ease in notation, we introduce
\begin{equation*}
    \|\sigma'\|_\infty:=\sup_{x\in\R^d}\|W\,D\sigma(x)\|,\quad\Gamma:=-\lambda_{\min}(D)+\|\sigma'\|_\infty,\quad\Gamma_1:=\lambda_{\max}(D)+\|\sigma'\|_\infty,\quad\Gamma_2:=\sup_{x\in\R^d}\|W\,D^2\sigma(x)\|. 
\end{equation*} 
Then, one has the following result, the proof of which is given in Section~\ref{ss:proof of Lipschitz constant of N_i for HNN}.
\begin{proposition}\label{pro:Lipschitz constant of N_i for HNN}
    Let $x^0\in\R^d$ and fix $i\in\{1,2\}$. Then $\cN_{\tau_i}:L^\infty((0, T); \R^k)\to\cS_d^+(\R)$ 
    is globally Lipschitz, 
    \begin{equation}\label{eq:Lipschitz constant of N_1 for HNN}
        \|\cN_{\tau_1}(u)-\cN_{\tau_1}(v)\|\le\frac{\Gamma_2\,\|B\|_\infty^3}{6}\left(3e^{\Gamma_1T}+1\right)\left(\frac{e^{\Gamma_1T}-1}{\Gamma_1}\right)^3\,\|u-v\|_\infty,\qquad\,u,\,v\in\,L^\infty((0, T); \R^k). 
    \end{equation}
    \begin{equation}\label{eq:Lipschitz constant of N_2 for HNN}
        \|\cN_{\tau_2}(u)-\cN_{\tau_2}(v)\|\le\begin{cases}
            \displaystyle\frac{\Gamma_2\,\|B\|_\infty^3}{3}\left(\frac{e^{\Gamma\,T}-1}{\Gamma}\right)^3\,\|u-v\|_\infty&\quad\text{if}\quad\Gamma\neq\,0,\cr
            \cr
            \displaystyle\frac{\Gamma_2\,\|B\|_\infty^3\,T^3}{3}\,\|u-v\|_\infty&\quad\text{if}\quad\Gamma=\,0,
        \end{cases}\qquad\,u,\,v\in\,L^\infty((0, T); \R^k).
    \end{equation}
\end{proposition}
\begin{remark}
    The bounds related to $\cN_{\tau_2}$ may involve a rate $\Gamma\le 0$, and are therefore sharper than the general estimates of Remark~\ref{rmk:specific case of Lipschitz constants}, where a direct application necessarily involves $\Gamma_1>0$ as the exponential rate.
\end{remark}
\begin{example}
    Consider the $2$D Hopfield-type recurrent neural network
    \begin{equation}\label{eq:example HRNN}
        \dot{x}=-D\,x+W\,\sigma(x)+Bu,\qquad x(0)=x^0\in\R^2
    \end{equation}
    where $x=(x_1,x_2)^\top\in\R^2$, $u\in L^\infty((0,T);\R)$, $B=(1, 0)^\top$, $D=\operatorname{diag}(d_1,d_2)$ with $d_i>0$, $W=(w_{ij})_{1\le i,j\le 2}$, and $\sigma_i(s)=\tanh(s)$. By~\cite[Theorem 2.1]{sun2007necessary}, the necessary and sufficient condition for the complete controllability of~\eqref{eq:example HRNN} is that the criterion function
    \[
    C(\gamma(t)) = -d_2\gamma_2(t)+w_{21}\tanh(\gamma_1(t))+w_{22}\tanh(\gamma_2(t))
    \]
    changes its sign over any control curve $\gamma(t):=(\gamma_1(t),\gamma_2(t))$ solution of
    \[
    \dot{\gamma}=1,\quad\dot{\gamma}_2=0\quad\Longleftrightarrow\quad\gamma_1(t)=t+\gamma_1^0,\quad\gamma_2(t)=\gamma_2^0,\quad\,t\in\R,\quad(\gamma_1^0,\gamma_2^0)\in\R^2.
    \]
    Since $d_2>0$, and $\tanh$ is bounded, $C(\gamma(t))$ (as a function of $t\in\R$) does not change its sign for large $\gamma_2^0$. Therefore, \eqref{eq:example HRNN} is not completely controllable. It follows that, neither the HCM~\cite{chitour2006continuation,schmoderer2022study,sussmann1992new,zhengping2024regularized} nor~\cite[Theorem 2.1]{sun2007necessary} can be directly applied for the motion planning of~\eqref{eq:example HRNN}. However, one proves using the Kalman-rank condition that $\cN_{\tau_i}(u)$ is invertible for all $(x^0, u)\in\R^2\times L^\infty((0,T);\R)$ whenever $w_{21}\ne 0$. In particular, Corollary~\ref{cor:relative coercivity from a zero control reference} applies to~\eqref{eq:example HRNN} for all $(x^0, x^1)\in\R^2\times\R^2$ such that $y_2:=x^1-\phi_T(x^0)$ satisfies~\eqref{eq:sufficient conditions for invertibility zero control}.
\end{example}

\section{Drift–modulated control-affine dynamics with time-varying and state–dependent input matrix}\label{s:sufficient conditions for nonautonomous quadratic nonlinear control systems with a pertubation} 

In this section, we consider the more general control-affine system
\begin{equation}\label{eq:nonlinear control system general bis}
    \dot{x}(t) = A(t,x(t))N_t(x(t))+B(t,x(t))u(t),\quad x(t_0) = x^0,\quad t\in[t_0, T]. 
\end{equation}

Throughout the following, $(t,x)\mapsto N_t(x)$ satisfies Assumption~\ref{ass:general assumption nonautonomous}. Whereas, the matrices $A$ and $B$ satisfy the following relaxed regularity assumptions:
\begin{assumption}\label{ass:gneneral assumption on A and B}
$A:[t_0, T]\times\R^d\to\R^{d\times d}$ is an element of $L^\infty((t_0, T)\times\R^d; \R^{d\times d})$, the input matrix $B:[t_0, T]\times\R^d\to\R^{d\times k}$ is an element of $L^\infty((t_0, T)\times\R^d; \R^{d\times k})$. Additionally, we assume that $x\mapsto A(\cdot,x),\,B(\cdot,x)$ are continuous and locally Lipschitz.
\end{assumption}

Note that under Assumptions~\ref{ass:general assumption nonautonomous} and~\ref{ass:gneneral assumption on A and B}, the Cauchy-Lipschitz theory guarantees the existence of a unique absolutely continuous solution $x\in C^0([t_0, T]; \R^d)$ to~\eqref{eq:nonlinear control system general bis} for any control $u\in\,L^\infty((t_0, T); \R^k)$; see, for instance, \cite[Chapter 2]{bressan2007introduction}. Observe also that in contrast to Section~\ref{s:state-dependent-input}, here we have only mild smoothness assumptions on $A$ and $B$. Still, under the trajectory freezing and compactness argument, we illustrate how the control analysis and synthesis of the baseline system~\eqref{eq:state-dependent-input} developed in Section~\ref{ss:controllability results nonautonomous} can be leveraged to provide some insight into the controllability properties of the general control-affine system~\eqref{eq:nonlinear control system general bis}. 

Fix $z\in C^0([t_0, T]; \R^d)$ and define $N_t^z$ and $B_z$ as
\begin{equation}
    N_t^z(x) := A(t,z(t))N_t(x),\quad B_z(t) := B(t,z(t))\qquad\forall (t,x)\in[t_0, T]\times\R^d.
\end{equation}
Then, $(t,x)\mapsto N_t^z(x)$ satisfies the same regularity assumptions as $(t,x)\mapsto N_t(x)$, as stated in Assumption~\ref{ass:general assumption nonautonomous} while $B_z$ satisfies Assumption~\ref{ass:gneneral assumption on B}. Additionally, the following estimates hold 
\begin{equation}
\|DN_t^z(w)\|\le\Lambda_3:=\Lambda_1\|A\|_\infty,\quad\|D^2N_t^z(w)\|\le\Lambda_4:=\Lambda_2\|A\|_\infty\quad\forall(t,w,z)\in[t_0,T]\times\R^d\times\,C^0([t_0, T]; \R^d)
\end{equation}
showing that $w\mapsto N_t^z(w)$ is globally $\Lambda_3$-Lipschitz continuous on $\R^d$, uniformly with respect to $(t,z)$. 

If $t\in[t_0, T]\to \Phi_{t_0,t}^z$ denotes the nonautonomous flow associated with $N_t^z$, then $\{\Phi_{s,t}^z\mid (s,t)\in[t_0, T]\}$ is a two-parameter family of diffeomorphisms which satisfies the same properties as the flow $\Phi_{t_0,t}$ of $N_t$ recalled in Section~\ref{s:prerequisites}. In particular, we have a similar lemma related to $\Phi_{s, t}^z$ as that given in Lemma~\ref{lem:spectral norm of D_s_t_phi and D^2_s_t_phi}.

Consider the following control-affine system  
\begin{equation}\tag{\textrm{$\Sigma_z$}}\label{eq:nonlinear control system general bis intermediate}
    \dot{x}(t) = N_t^z(x(t))+B_z(t)u(t),\quad x(t_0) = x^0,\quad t\in[t_0, T]. 
\end{equation}

As for \eqref{eq:state-dependent-input}, system~\eqref{eq:nonlinear control system general bis intermediate} has unique absolutely continuous solution $x_u^z\in C^0([t_0, T]; \R^d)$ that can be represented as in Theorem~\ref{thm:sol representation nonautonomous} with $b(t,x) = B_z(t)u(t)$. We use the backward representation and get
\begin{equation}
    x_u^z(t)=\Phi_{T,t}^z\!\left(\Phi_{t_0,T}^z(x^0)+\int_{t_0}^{t} D\Phi_{s,T}^z\big(x_u^z(s)\big)\,B_z(s)\,u(s)\,ds\right),
\qquad t\in[t_0,T]. 
\end{equation}

Fix $x^0\in\R^d$, $u\in L^\infty((t_0, T); \R^k)$ and let $x_u^z:=x_{u,x^0}^z\in C^0([t_0, T]; \R^d)$ be the corresponding solution of  \eqref{eq:nonlinear control system general bis intermediate}. We define $\cN_{\tau_2}^z(u):=\cN_{\tau_2}^z(u, x^0)$ similarly as was defined $\cN_{\tau_2}(u)$ in Section~\ref{ss:controllability results nonautonomous}, viz.
\begin{equation}
    \cN_{\tau_2}^z(u) = \int_{t_0}^{T}D\Phi_{t,T}^z(x_u^z(t))B_z(t)B_z(t)^\top D\Phi_{t,T}^z(x_u^z(t))^\top\,dt,
\end{equation}
which is a symmetric positive-semidefinite matrix. In particular, Remark~\ref{rmk:specific case of Lipschitz constants} yields
 \begin{equation}\label{eq:Lipschitz constant of N_2^z}
        \|\cN_{\tau_2}^z(u)-\cN_{\tau_2}^z(v)\|\le\frac{\Lambda_4\,\|B\|_\infty^3}{3}\left(\frac{e^{\Lambda_3\Delta t_0}-1}{\Lambda_3}\right)^3\,\|u-v\|_\infty\qquad\forall\,u,\,v\in\,L^\infty((t_0,T); \R^k),
    \end{equation}
showing that $\cN_{\tau_2}^z:L^\infty((t_0, T); \R^k)\to\cS_d^+(\R)$ is uniformly globally Lipschitz w.r.t. $z\in C^0([t_0, T]; \R^d)$.

Similarly, as for
$W_2(T)$ defined in \eqref{eq:SSD W_i}, one introduces the matrix 
\begin{equation}\label{eq:SSD W_iz}
W_2^z(T):=\cN_{\tau_2}^z(0)=\int_{t_0}^TD\Phi_{t,T}^z(\Phi_{t_0,t}^z(x^0))B_z(t)B_z(t)^\top D\Phi_{t,T}^z(\Phi_{t_0,t}^z(x^0))^\top\,dt,
\end{equation}
which is symmetric positive semi-definite.

The main theorem of this section is then the following. 
\begin{theorem}\label{thm:main general}
    Let $x^0\in\R^d$. Assume that there exists $C>0$ such that
    \begin{equation}\label{eq:boundeness of lambda_min}
        \lambda_{\min}(W_2^z(T))\ge C^{-1}\qquad\forall z\in C^0([t_0, T]; \R^d),
    \end{equation}
    and that, for each frozen system, the corresponding synthesis map satisfies the contractivity condition of Theorem~\ref{thm:abstract-fixed-point}.
    Then, there exist $z_{*}\in C^0([t_0, T]; \R^d)$ such that \eqref{eq:nonlinear control system general bis} is controllable on $[t_0, T]$ from $x^0$ to any target state $x^1\in\R^d$ satisfying for some $\theta\in(0, 1)$ the following estimate
    \begin{equation}\label{eq:estimate reachable set general}
        |x^1-\Phi_{t_0,T}^{z_*}(x^0)|\le\frac{\theta\,\lambda_{\min}\big(W_2^{z_*}(T)\big)}{(1+\theta)\,L_{*}\,\|B\|_\infty\,e^{\Lambda_3\Delta t_0}}.
    \end{equation}    
    Here $L_*>0$ is the Lipschitz constant in~\eqref{eq:Lipschitz constant of N_2^z}.
\end{theorem}

\begin{remark}
Following Theorem~\ref{thm:relative coercivity from a control reference}, the assumption in Theorem~\ref{thm:main general} can be relaxed as follows: assume the existence of a reference control $u_2^\flat\in L^\infty((t_0,T);\R^k)$ and a constant $C>0$ such that
\[
\lambda_{\min}\big(\cN_{\tau_2}^z(u_2^\flat)\big)\ \ge\ C^{-1},\qquad \forall\, z\in C^0([t_0,T];\R^d).
\]
Moreover, as illustrated in Section~\ref{ss:refined estimates for HNN}, the constant $L_\ast>0$ is often conservative, and can be sharpened in concrete models by exploiting structural properties of the vector field $N_t^z$.
\end{remark}

\begin{proof}[\textit{Proof} of Theorem~\ref{thm:main general}]
Let $z\in C^0([t_0, T]; \R^d)$. By~\eqref{eq:boundeness of lambda_min}, $W_2^z(T)$ is invertible. Then, letting $C_2:=(1+\theta)/C$ with $C>0$ as in~\eqref{eq:boundeness of lambda_min}, and for some $\theta\in(0,1]$, the feasible coercivity class
\begin{equation}\label{eq:ACC z}
    \cF(C_2) = \{u\in L^\infty((t_0, T); \R^k):\lambda_{\min}(\cN_{\tau_2}^z(u))\ge C_2^{-1}\}
\end{equation}
 is nonempty since it contains $u=0$. For $x^1\in\R^d$, if the following estimate holds
\begin{equation}
    |x^1-\Phi_{t_0,T}^{z}(x^0)|\le\frac{\theta\,\lambda_{\min}\big(W_2^{z}(T)\big)}{(1+\theta)\,L_{*}\,\|B\|_\infty\,e^{\Lambda_3\Delta t_0}},
\end{equation} 
then by Remark~\ref{rmk:specific case of Lipschitz constants}, and Theorems~\ref{thm:abstract-fixed-point} and~\ref{thm:main controllability result nonlinear}, the fixed point $u_z\in\cF(C_2)$ of
\[
\cS_2^z(u)(t)=B_z(t)^\top D\Phi_{t,T}^z(x_{u}(t))^\top(\cN_{\tau_2}^z(u))^{-1}\left(x^1-\Phi_{t_0,T}^z(x^0)\right)
\]
exists, is unique, and the corresponding solution $x_{u_z}$ to system~\eqref{eq:nonlinear control system general bis intermediate} satisfies $x_{u_z}(T)=x^1$.

To complete the proof of the theorem, let us define the map
    \begin{equation}\label{eq:map cF}
        \cZ:z\in C^0([t_0, T]; \R^d)\mapsto\cZ(z)=x_{u_z}\in C^0([t_0, T]; \R^d)
    \end{equation}
    where $x_{u_z}\in C^0([t_0, T]; \R^d)$ is the solution of \eqref{eq:nonlinear control system general bis intermediate} corresponding to the control $u_z=\cS_2^z(u_z)$. 
    
If $\cZ$ has a fixed point $z_*\in C^0([t_0, T]; \R^d)$, then $z_*=x_{u_{z_*}}$, and by construction, $u_{z_*}$ steers \eqref{eq:nonlinear control system general bis} from $x^0$ to any target state $x^1\in\R^d$ satisfying~\eqref{eq:estimate reachable set general}, i.e., $x_{u_{z_*}}$ is the solution to system~\eqref{eq:nonlinear control system general bis} and satisfies $x_{u_{z_*}}(T)=x^1$.
    
First, $\cZ$ is clearly well-defined and continuous. 
Moreover, by Lemma~\ref{lem:spectral norm of D_s_t_phi and D^2_s_t_phi} and \eqref{eq:boundeness of lambda_min}, it holds that
\begin{equation}\label{eq:boundeness of u_z}
    \|u_z\|_\infty\le C\|B\|_\infty e^{\Lambda_3\Delta t_0}(e^{\Lambda_3\Delta t_0}|x^1|+|x^0|)\quad\forall z\in C^0([t_0, T];\R^d) 
\end{equation}
so that for some $M>0$ independent of $z$, we find
\begin{equation}\label{eq:boundeness of cF}
    \|\cZ(z)\|_\infty=\|x_{u_z}\|_\infty\le e^{\Lambda_1\Delta t_0}(|x^0|+\Delta t_0(\|B\|_\infty\|u_z\|_\infty))\le M\quad\forall z\in C^0([t_0, T];\R^d). 
\end{equation}
Next, $\cZ$ is compact. In fact, let $\cB\subset C^0([t_0, T]; \R^d)$ be a bounded set, let us show that $\overline{\cZ(\cB)}$ is compact. First, $\cZ(\cB)\subset C^0([t_0, T]; \R^d)$ is bounded by \eqref{eq:boundeness of cF}. It remains to prove that for all $t_*\in [t_0, T]$, the following holds
\begin{equation}\label{eq:uniform equicontinuity z}
    \sup\limits_{z\in\cB}|\cZ(z)(t)-\cZ(z)(t_*)|\xrightarrow[t\to t_*]{} 0.
\end{equation}
One has immediately from \eqref{eq:boundeness of u_z} and \eqref{eq:boundeness of cF} that for some $M_1>0$, it holds that
\begin{equation}\label{eq:norm of Y(u)(t)-Y(u)(t_0) z}
     |\cZ(z)(t)-\cZ(z)(t_*)|=|x_{u_z}(t)-x_{u_z}(t_*)|\le(\Lambda_3\|x_{u_z}\|_\infty+\|B\|_\infty\|u_z\|_\infty)|t-t_*|\le M_1|t-t_*| 
\end{equation}
for every $t\in[t_0, T]$. This proves \eqref{eq:uniform equicontinuity z} and, therefore, that $\cZ$ is compact by the Ascoli-Arzelà theorem. Finally, the set $\{z\in C^0([t_0, T]; \R^d)\mid\, z=\theta \cZ(z)\,\text{for some}\, \theta\in [0, 1]\}$ is bounded by~\eqref{eq:boundeness of u_z} and~\eqref{eq:boundeness of cF}. According to Schaefer’s fixed point theorem, $\cZ$ admits at least one fixed point $z_*\in C^0([t_0, T]; \R^d)$, i.e., $z_* = \cF(z_*) = x_{u_{z_*}}$.
\end{proof}

\begin{remark}
Theorem~\ref{thm:main general} generalizes~\cite[Theorem~3.40]{coron2007control} beyond the trivial case $N_t(x)=x$.  
In this setting, $N_t^z(x)=A(t,z(t))x$, so $\cN_{\tau_2}^z(u)$ is independent of $u$ and the admissible class in~\eqref{eq:ACC z} reduces to $\cF(C_2)=L^\infty((t_0,T);\R^k)$ for all initial states and trajectories by letting $C_2:=C$ with $C>0$ as in~\eqref{eq:boundeness of lambda_min}.  
Hence, by arguments analogous to Proposition~\ref{pro:uniformly invertible B}, the baseline system~\eqref{eq:nonlinear control system general bis intermediate} is globally controllable for all $z$, and global controllability of~\eqref{eq:nonlinear control system general bis} follows via Schaefer’s fixed point theorem.
\end{remark}

\section{Conclusion}\label{s:conclusion}
We developed a constructive framework of control synthesis for control–affine system $\dot x=N_t(x)+B(t,x)u$ based on trajectory–dependent Gramians and fixed–point synthesis maps. On a ball $\cF(y_i)$ of the feasible coercivity class $\cF(C_i)$, sufficient conditions ensure that the synthesis map is a self–map and admits a unique fixed point with the associated control satisfying an explicit energy certificate. Beyond synthesis, we established existence of minimizers: the feasible set $\mathfrak F_i=\{u\in\cF(y_i):\,L_{u,\tau_i}u=y_i\}$ is weakly sequentially closed in $L^2((t_0, T); \R^k)$, hence $\tfrac12\|u\|_{L^2}^2$ attains a minimum on $\mathfrak F_i$ (Proposition~\ref{prop:exist-on-ball-global}). We then linked optimality to fixed points: under the orthogonality condition at the fixed point (Theorem~\ref{thm:minimizer-is-fp}) and a finite-dimensional invertibility assumption, every local minimizer over $\mathfrak F_i$ equals the fixed point $u_i$; therefore, the minimizer on $\mathfrak F_i$ is unique and global. If, in addition, $\cF(C_i)=L^\infty((t_0, T);\R^k)$, the same conclusion holds among all bounded controls satisfying the endpoint constraint.

The framework extends the classical Gramian theory for linear systems and, for~\eqref{eq:nonlinear control system general}, recovers~\cite[Theorem.~3.40]{coron2007control} in the case of $N_t(x)=x$. Moreover, our results apply to Hopfield–type networks~\cite{hopfield1984neurons} and Lur’e systems~\cite{diwadkar2015control} in general. It is constructive and amenable to numerical implementation. The analysis is pointwise in the initial state; natural next steps include deriving uniform state–space conditions paralleling the linear equivalence between controllability and Gramian invertibility, investigating robustness, and obtaining a sharper quantitative characterization of the reachable set.

When complete controllability on $[t_0, T]$ is known \emph{a priori} for system $ \dot{x}=N_t(x)+B(t, x)u$, our framework suggests a systematic procedure: certify on short windows that (i) $\cF(C_i)$ is nonempty (via an invertible reference Gramian), and (ii) the synthesis map is a self–map and satisfies an iterate-contractivity condition; then concatenate the windowed fixed points along a partition of $[t_0, T]$ to steer arbitrary targets, with additive energy bounds. This ``windowed Gramian synthesis'' provides an implementable alternative to homotopy continuation, equipped with an energy certificate and an \emph{a posteriori} optimality upgrade at the concatenated fixed point via the orthogonality condition test. Establishing uniform window conditions from structural data and extending the patching argument to underactuated regimes are promising directions for future work.


\appendix

\section{Proofs of results stated in the main text}\label{app:proofs of some results presented in the main text}
In this section, we present proofs of some of the results from the main text. We start with the following.

\subsection{Proof of Theorem~\ref{thm:sol representation nonautonomous}}\label{ss:proof of sol representation nonautonomous}
We present in this section the proof of the solution representation \eqref{eq:forward-backward representation}.

\begin{proof}
    Under the assumptions of Theorem~\ref{thm:sol representation nonautonomous}, system 
    \begin{equation}\label{eq:local}
        \dot{x}(t)=N_t(x(t))+b(t, x(t)),\qquad x(t_0)=x^0    
    \end{equation}
    has (see, for instance, \cite[Chapter 2]{bressan2007introduction}) a unique absolutely continuous solution $x\in C^0([t_0, T]; \R^d)$.
Let us prove that this solution can be represented by~\eqref{eq:forward-backward representation}. Let $t\in[t_0, T]\mapsto y(t)\in\R^d$ be such that $t\mapsto x(t) = \Phi_{\tau_i,t}(y(t))$ is the solution of \eqref{eq:local}. Then $y$ is derivable almost everywhere w.r.t. $t$ and it holds that
    \begin{equation}
        N_t(x(t))+b(t, x(t))=\dot{x}(t)=N_t(x(t))+D\Phi_{\tau_i,t}(y(t))\,\dot{y}(t)
    \end{equation}
    so that using \eqref{eq:inverse Dphi_t_0,t}, we find that $y$ solves the following.
    \begin{equation}\label{eq:equation satisfying b}
        \dot{y}(t)=D\Phi_{t,\tau_i}(x(t))\,b(t, x(t)),\qquad y(t_0)=\Phi_{t_0,\tau_i}(x^0).
    \end{equation}
    Integrating \eqref{eq:equation satisfying b} yields \eqref{eq:forward-backward representation}. Conversely, if \eqref{eq:forward-backward representation} holds, then $x(t_0)=\Phi_{\tau_i,t_0}(\Phi_{t_0,\tau_i}(x^0))=x^0$ and $x\in C^0([t_0, T]; \R^d)$ by composition. Otherwise, there exists $(t_n)\subset [t_0, T]$, $t_{*}\in[t_0, T]$ with $t_n\to t_{*}$ and $\varepsilon>0$ such that $|\Phi_{t,\tau_i}(x(t_n))-\Phi_{t_{*},\tau_i}(x(t_{*}))|\ge\varepsilon$. In fact, $x\in C^0([t_0, T]; \R^d)$ if and only if $t\mapsto y(t):=\Phi_{t,\tau_i}(x(t))$ belongs to $C^0([t_0, T]; \R^d)$ since $\Phi_{t,\tau_i}$ is invertible and $C^1$ w.r.t. $t\in\R$. However, we have from \eqref{eq:forward-backward representation} that
    \[
    \Phi_{t_n,\tau_i}(x(t_n))-\Phi_{t_{*},\tau_i}(x(t_{*})) = \int_{t_*}^{t_n}D\Phi_{s,\tau_i}(x(s))\,b(s,x(s))\,ds
    \]
    which implies $|\Phi_{t_n,\tau_i}(x(t_n))-\Phi_{t_{*},\tau_i}(x(t_{*}))|\le \|b\|_\infty e^{\Lambda_1|\tau_i-t_*|}|t_n-t_*|$
    by using \eqref{eq:EGectral norm of D_s_t_phi} with $\beta(s)=x(s)$. It follows that $|\Phi_{t_n,\tau_i}(x_u(t_n))-\Phi_{t_{*},\tau_i}(x_u(t_{*}))|\to 0$ as $n\to\infty$, which is inconsistent. 
    Letting now
    \begin{equation}
        z(t):= \Phi_{t_0,\tau_i}(x^0)+\int_{t_0}^{t}D\Phi_{s,\tau_i}(x(s))\,b(s,x(s))\,ds
    \end{equation}
    and deriving \eqref{eq:forward-backward representation} almost everywhere w.r.t. $t$ yields
    \begin{equation}\label{eq:derivation a}
        \dot{x}(t)=\partial_t\Phi_{\tau_i,t}(z(t))+D\Phi_{\tau_i,t}(z(t))\,D\Phi_{t,\tau_i}(x(t))\,b(t,x(t))=N_t(\Phi_{\tau_i,t}(z(t)))+b(t,x(t))=N_t(x(t))+b(t,x(t))
    \end{equation}
   by $x(t) = \Phi_{\tau_i,t}(z(t))$ and $D\Phi_{\tau_i,t}(z(t))\,D\Phi_{t,\tau_i}(x(t))=\idty$. It follows that \eqref{eq:forward-backward representation} solves \eqref{eq:local}.
\end{proof}

\subsection{Proof of Lemma~\ref{lem:state-transition matrix of the nonlinear control B=0}}\label{ss:state-transition matrix of the nonlinear control}
\begin{proof}
Recall from \cite[Theorem~3.2.6]{bressan2007introduction} that the Fréchet derivative $D_ux_u(t)h$ of $x_u(t)$ w.r.t. $u$ in the direction of $h\in L^\infty((t_0, T); \R^k)$ is given by
\begin{equation}\label{eq:F-derivative of x_u}
    D_ux_u(t)h = \int_{t_0}^t\!\!\!R_u(t,s)B(s,x_u(s))h(s)\,ds.
\end{equation}
Using the solution representation~\eqref{eq:forward-backward representation} with $b(t,x(t))=B(t,x(t))u(t)$
and letting $y_u(t):=\Phi_{t,\tau_i}(x_u(t))$, one finds $D_ux_u(t)h=\big[D\Phi_{t,\tau_i}(x_u(t))\big]^{-1}D_uy_u(t)h$ and
\[
\frac{d}{dt}D_uy_u(t)h = \cA_{\tau_i,u}(t)D_uy_u(t)h+D\Phi_{t,\tau_i}(x_u(t))B(t,x_u(t))h(t),\qquad D_uy_u(t)h|_{t=t_0}=0
\]
where $\cA_{\tau_i,u}(t)$ is defined as in the statement of Lemma~\ref{lem:state-transition matrix of the nonlinear control B=0}.
It follows that
\begin{equation}
    D_uy_u(t)h = \int_{t_0}^t\!\!M_{\tau_i,u}(t,s)D\Phi_{s,\tau_i}(x_u(s))B(s,x_u(s))h(s)\,ds
\end{equation}
where $M_{\tau_i,u}(t,s)$ is defined as in the statement of Lemma~\ref{lem:state-transition matrix of the nonlinear control B=0}. One deduces that
\begin{equation}\label{eq:F-derivative of x_u bis}
    D_ux_u(t)h = \big[D\Phi_{t,\tau_i}(x_u(t))\big]^{-1}\int_{t_0}^t\!\!M_{\tau_i,u}(t,s)D\Phi_{s,\tau_i}(x_u(s))B(s,x_u(s))h(s)\,ds.
\end{equation}
Next, \eqref{eq:F-derivative of x_u} and~\eqref{eq:F-derivative of x_u bis} suggest the factorization $R_u(t,s) = \big[D\Phi_{t,\tau_i}(x_u(t))\big]^{-1}M_{\tau_i,u}(t,s)D\Phi_{s,\tau_i}(x_u(s))$. Using this factorization, one checks immediately that
\[
R_u(t_1,t_1)=\idty,\quad 
R_u(t_1,t_2)R_u(t_2,t_3)=R_u(t_1,t_3),\quad
R_u(t_1,t_2)R_u(t_2,t_1)=\idty,\qquad \forall(t_1,t_2,t_3)\in[t_0,T]^3
\]
since $M_{\tau_i,u}(t,s)$ is a state-transition matrix. Finally, using the following identities
\begin{eqnarray}
\frac{d}{dt}\big[D\Phi_{t,\tau_i}(x_u(t))\big]^{-1}&=&-\big[D\Phi_{t,\tau_i}(x_u(t))\big]^{-1}\left(\frac{d}{dt}D\Phi_{t,\tau_i}(x_u(t))\right)\big[D\Phi_{t,\tau_i}(x_u(t))\big]^{-1}\nonumber\\
&=&DN_t(x_u(t))\big[D\Phi_{t,\tau_i}(x_u(t))\big]^{-1}+D^2\Phi_{\tau_i,t}\big(\Phi_{t,\tau_i}(x_u(t))\big)D\Phi_{t,\tau_i}(x_u(t))B(t,x_u(t))u(t),
\end{eqnarray}
\begin{equation}
    D^2\Phi_{\tau_i,t}\big(\Phi_{t,\tau_i}(x_u(t))\big)D\Phi_{t,\tau_i}(x_u(t))D\Phi_{t,\tau_i}(x_u(t))+ D\Phi_{\tau_i,t}\big(\Phi_{t,\tau_i}(x_u(t))\big)D^2\Phi_{t,\tau_i}(x_u(t))=0,
\end{equation}
and the fact that $M_{\tau_i,u}(t,s)$ is a state-transition matrix of~\eqref{eq:state-transition matrix M_T,u}, one finds that
\begin{equation}\label{eq:resolvent properties}
\frac{\partial R_u}{\partial t}(t,s)=\big[DN_t(x_u(t))+D_xB(t,x_u(t))u(t)\big]R_u(t,s),\qquad\quad (t,s)\in[t_0,T]^2.
\end{equation}
This completes the proof of the Lemma.
\end{proof}

\subsection{Proof of Theorem~\ref{thm:abstract-fixed-point}}\label{ss:abstract-fixed-point}
In this section, we present the proof of Theorem~\ref{thm:abstract-fixed-point} that we split into several steps for the reader's convenience. Since the maps $\cS_1$ and $\cS_2$ play similar roles, we only focus on $\cS_2$. Throughout this section, we set $\cY:=L^\infty((t_0, T); \R^k)$. Under the assumptions of Theorem~\ref{thm:abstract-fixed-point},
\[
\cS_2(u) = B_u(t)^\top\,Q_u(t)^\top\,\cN_{\tau_2}(u)^{-1}\,y_2\qquad\text{and}\qquad \cS_2(u)\in\cF(y_2)\qquad\forall u\in\cF(y_2),
\]
where $y_2:=x^1-\Phi_{t_0,T}(x^0)$, $B_u(t):=B(t,x_u(t))$ and $Q_u(t):=D\Phi_{t,T}(x_{u}(t))$. 
For ease in notation, for $s\in [0, \Delta t_0]$, where $\Delta t_0:=T-t_0$, we let
\begin{equation}\label{eq:cst1}
      E_0=e^{L_B\,\zeta_2\,\Delta t_0},\;\;E_1(s) = e^{\Lambda_1s},\;\, E_2(s)=e^{2\Lambda_1s}-e^{\Lambda_1s},\;\, E_{1,\infty}=\max_{s\in[0, \Delta t_0]}E_1(s),\;\, E_{2,\infty}=\max_{s\in[0, \Delta t_0]}E_2(s).
\end{equation}

We introduce the following Volterra and Fredholm-type bounded linear operators $V,\;F:\cY\to \cY$ defined by
\begin{equation}
    (Vf)(t) = \int_{t_0}^tf(\tau)\,d\tau,\qquad Ff = \int_{t_0}^Tf(\tau)\,d\tau\qquad\forall (t,\,f)\in[t_0, T]\times\cY.
\end{equation}
\begin{lemma}\label{lem:state-and-gramian-estimates-corrected}
    For every $(u,\,v)\in\cF(y_2)^2$, the following estimates hold
    \begin{equation}\label{eq:estimate x_u-x_v}
        |x_u(t)-x_v(t)|\le\|B\|_\infty\, E_0\,E_{1,\infty}\,(V|u-v|)(t)\le \|B\|_\infty\,E_0\,E_{1,\infty}\,\Delta t_0\|u-v\|_\infty\qquad\forall t\in [t_0, T],
    \end{equation}
    \vspace{-0.4cm}
    \begin{equation}\label{eq:estimate N_2(u)-N_2(v)}
        \|\cN_{\tau_2}(u)-\cN_{\tau_2}(v)\|\le2\,E_0\,E_{1,\infty}\,\|B\|_\infty^2\bigg(\frac{\|B\|_\infty\Lambda_2\,E_{2,\infty}}{\Lambda_1}+L_B\,E_{1,\infty}\bigg)\,F\circ V|u-v|. 
    \end{equation}
\end{lemma}
\begin{proof}
    Recall that $\dot{x}_u(t)=N_t(x_u(t))+B(t,x_u(t))u(t)$, $x_u(t_0)=x^0$ and $\dot{x}_v(t)=N_t(x_v(t))+B(t,x_u(t))v(t)$, $x_v(t_0)=x^0$. Then, by the Cauchy-Schwarz inequality, one finds that
    \[
    \frac{d}{dt}|x_u(t)-x_v(t)|=\frac{\langle\dot{x}_u(t)-\dot{x}_v(t),x_u(t)-x_v(t)\rangle}{|x_u(t)-x_v(t)|}\le(\Lambda_1+L_B\,\zeta_2)\,|x_u(t)-x_v(t)|+\|B\|_\infty\,|u(t)-v(t)|
    \]
    which by Gronwall's lemma yields~\eqref{eq:estimate x_u-x_v}. Now using Lemma~\ref{lem:spectral norm of D_s_t_phi and D^2_s_t_phi} respectively with $\beta\in\{x_u,x_v\}$, one finds
    \begin{equation}
        \|\cN_{\tau_2}(u)-\cN_{\tau_2}(v)\|\le\,2\bigg(\frac{\|B\|_\infty^2\Lambda_2\,E_{2,\infty}}{\Lambda_1}+L_B\|B\|_\infty\,E_{1,\infty}\bigg)\int_{t_0}^T\,|x_u(s)-x_v(s)|\,ds
    \end{equation}
    which by the first inequality in~\eqref{eq:estimate x_u-x_v} yields~\eqref{eq:estimate N_2(u)-N_2(v)}. 
\end{proof}

We introduce the positive constants
\begin{equation}\label{eq:alpha 1 and 2}
\alpha_1=2E_0C_2^2\|B\|_\infty^3E_{1,\infty}^2\bigg(\frac{\Lambda_2\|B\|_\infty E_{2,\infty}}{\Lambda_1}+L_BE_{1,\infty}\bigg)|y_2|,\;\alpha_2=E_0C_2\|B\|_\infty E_{1,\infty}\bigg(\frac{\Lambda_2\|B\|_\infty E_{2,\infty}}{\Lambda_1}+L_BE_{1,\infty}\bigg)|y_2|,
\end{equation}
and the bounded and linear operator $\cL:\cY\to\cY$, defined by
\begin{equation}\label{eq:cL}
    (\cL f)(t) = \alpha_1\,F\circ V f+\alpha_2(V f)(t)\qquad\forall (t,\,f)\in[t_0, T]\times\cY.
\end{equation}

\begin{lemma}\label{lem:norm-of-cL-corrected}
Set
\begin{equation}\label{eq:K2-p2}
    K_2:=\Delta t_0(\alpha_1\Delta t_0+\alpha_2),\qquad
    p_2:=\frac{\alpha_1\Delta t_0}{\alpha_1\Delta t_0+\alpha_2}\in(0,1).
\end{equation}
Let $\gamma_m(p_2):=z_m(1)$, where $z_0\equiv1$ on $[0,1]$ and
\begin{equation}\label{eq:gamma-recursion}
    z_{m+1}(y)=(1-p_2)\int_0^y z_m(r)\,dr
    +p_2\int_0^1(1-r)z_m(r)\,dr.
\end{equation}
Then, for every $m\ge1$,
\begin{equation}\label{eq:norm-of-cL-corrected}
    \|\cL^m\|\le K_2^m\gamma_m(p_2).
\end{equation}
\end{lemma}
\begin{proof}
Introduce the isometric rescaling $\mathcal Rf(y):=f(t_0+\Delta t_0y)$ on $L^\infty(0,1)$, and let $(Wg)(y)=\int_0^y g(r)dr$ and $Jg=\int_0^1g(r)dr$. Then $\mathcal RV\mathcal R^{-1}=\Delta t_0W$ and $\mathcal RF\mathcal R^{-1}=\Delta t_0J$. Hence $\widetilde{\cL}:=\mathcal R\cL\mathcal R^{-1}=K_2T_{p_2}$, with $T_{p_2}:=(1-p_2)W+p_2JW$.

Since $T_{p_2}$ is positive, for $\|g\|_\infty\le1$ one has $|T_{p_2}^mg|\le T_{p_2}^m\mathbf1$ a.e. Let $z_m:=T_{p_2}^m\mathbf1$. Then $z_0\equiv1$ and $z_{m+1}=T_{p_2}z_m$, which gives exactly~\eqref{eq:gamma-recursion}. The functions $z_m$ are nonnegative and nondecreasing by induction; hence $\|z_m\|_\infty=z_m(1)=\gamma_m(p_2)$. Thus $\|T_{p_2}^m\|\le\gamma_m(p_2)$, and since $\mathcal R$ is an isometry, $\|\cL^m\|=\|\widetilde{\cL}^{\,m}\|\le K_2^m\gamma_m(p_2)$.
\end{proof}

\begin{lemma}\label{lem:gamma-threshold}
Let $p\in(0,1)$, and let $\lambda(p)>0$ be the unique positive solution of
\begin{equation}\label{eq:lambda-p}
    p(e^\lambda-1)-\lambda=0.
\end{equation}
Set $\Lambda(p):=\lambda(p)/(1-p)$. Then $\Lambda(p)>2$. Moreover, if $0<K<\Lambda(p)$, then $K^m\gamma_m(p)\to0$ as $m\to\infty$.
\end{lemma}
\begin{proof}
For $T_p=(1-p)W+pJW$, the same argument as in the proof of Lemma~\ref{lem:norm-of-cL-corrected} gives $\|T_p^m\|\le\gamma_m(p)$, while testing on $\mathbf 1$ gives the reverse inequality. Hence $\gamma_m(p)=\|T_p^m\|$. The compact operator $T_p$ is strongly positive on the cone of nonnegative functions in $L^\infty(0,1)$, so the Krein--Rutman theorem~\cite[Theorem 6.13 \& Problem 41]{brezis2011functional} gives an eigenpair $\mu=r(T_p)>0$ and $f\in\operatorname{Int}P$ such that $T_pf=\mu f$. It follows that $f$ has an absolutely continuous representation on $[0,1]$ that we still denote by $f$. Differentiating $T_pf=\mu f$ gives $\mu f'=(1-p)f$, hence $f(y)=f(0)e^{\lambda y}$ with $\lambda=(1-p)/\mu$. Evaluating at $y=0$ gives~\eqref{eq:lambda-p}; hence $r(T_p)=1/\Lambda(p)$. The inequality $\Lambda(p)>2$ follows by setting $q=1-p$ and observing that the defining function is negative at $\lambda=2q$, equivalently $2q<\log((1+q)/(1-q))$. The spectral-radius formula gives $\lim_{m\to\infty}\gamma_m(p)^{1/m}=1/\Lambda(p)$, and therefore $K^m\gamma_m(p)\to0$ whenever $K<\Lambda(p)$.
\end{proof}

\begin{proof}[\textit{Proof} of Theorem~\ref{thm:abstract-fixed-point}]
    Let $(u,\,v)\in\cF(y_2)^2$. Then, for every $t\in [t_0, T]$, it hols that 
    \begin{eqnarray}
        |\cS_2(u)(t)-\cS_2(v)(t)|&\le&|y_2|C_2\big(\|B\|_\infty\,E_{1,\infty}\,C_2\|\cN_{\tau_2}(u)-\cN_{\tau_2}(v)\|+\|B\|_\infty\|Q_u(t)-Q_v(t)\|+E_{1,\infty}\|B_u(t)-B_v(t)\|\big)\nonumber\\
        &\le&\alpha_1\,F\circ V |u-v|+\alpha_2(V |u-v|)(t)=(\cL |u-v|)(t)
    \end{eqnarray}
    by Lemma~\ref{lem:spectral norm of D_s_t_phi and D^2_s_t_phi}, \eqref{eq:estimate x_u-x_v}, \eqref{eq:estimate N_2(u)-N_2(v)}, \eqref{eq:alpha 1 and 2} and~\eqref{eq:cL}. It follows that for every $m\in\N$, with $m\ge 1$, it holds that
    \begin{equation}
         |\cS_2^m(u)(t)-\cS_2^m(v)(t)|\le (\cL^m |u-v|)(t)\le\|(\cL^m |u-v|)\|_\infty\le\|\cL^m\|\|u-v\|_\infty
        \le K_2^m\gamma_m(p_2)\|u-v\|_\infty,
    \end{equation}
   for all $t\in[t_0, T]$, by Lemma~\ref{lem:norm-of-cL-corrected}. This proves~\eqref{eq:supnorm of F_i^m(u)-F_i^m(v)} with $\varrho_m:=K_2^m\gamma_m(p_2)$. If there exists $m_i\ge1$ such that $\varrho_{m_i}<1$, one completes the proof of the theorem by the Banach fixed-point theorem~\cite[Theorem~2.4]{almezel2014topics}.
\end{proof}

\subsection{Proof of Lemma~\ref{lem:preparatory result}}\label{ss:preparatory result}

For clarity, we split the proof into several steps. Throughout this section, we set $\cX:=L^2((t_0,T);\R^k)$, $\cY:=L^\infty((t_0,T);\R^k)$,
and let $\cF(C_i)$ be defined by~\eqref{eq:feasible coercivity class}. Fix
$u\in\cY$, and consider the following linear and bounded operators
$H_{u,\tau_i}$, $W_{u,\tau_i}$, and $J_{u,\tau_i}$ defined by
\begin{equation*}
    (H_{u,\tau_i}h)(t)
    =
    R_u(\tau_i,t)B(t,x_u(t))h(t),
    \qquad
    (W_{u,\tau_i}h)(t)
    =
    R_u(t,\tau_i)\int_{t_0}^t h(s)\,ds,
\end{equation*}
and
\begin{equation*}
    J_{u,\tau_i}h
    =
    \int_{t_0}^T
    D^2\Phi_{t,\tau_i}(x_u(t))B(t,x_u(t))u(t)h(t)\,dt
    +
    \int_{t_0}^T
    D\Phi_{t,\tau_i}(x_u(t))D_xB(t,x_u(t))u(t)h(t)\,dt .
\end{equation*}
Here $R_u(t,s)$ is defined in
Lemma~\ref{lem:state-transition matrix of the nonlinear control B=0}. Since
$R_u$ is a state transition matrix, one has
\[
(W_{u,\tau_i}H_{u,\tau_i}h)(t)
=
R_u(t,\tau_i)\int_{t_0}^t(H_{u,\tau_i}h)(s)\,ds
=
\int_{t_0}^tR_u(t,s)B(s,x_u(s))h(s)\,ds.
\]
We set $K_{u,\tau_i}:=J_{u,\tau_i}W_{u,\tau_i}H_{u,\tau_i}$. It is clear by Assumptions~\ref{ass:general assumption nonautonomous}
and~\ref{ass:gneneral assumption on B} that
$\cE\in C^1(\cY;\R^d)$; see, for instance,
\cite[Theorem~3.2.6]{bressan2007introduction}. Hence
$G_{\tau_i}\in C^1(\cY;\R^d)$ by composition. Moreover, for any
$h\in\cX$,
\begin{equation}\label{eq:F-derivatives}
    D\cE(u)h
    =
    \int_{t_0}^T
    R_u(T,t)B(t,x_u(t))h(t)\,dt,
    \qquad
    DG_{\tau_i}(u)h
    =
    L_{u,\tau_i}h+K_{u,\tau_i}h.
\end{equation}
Using the forward--backward representation~\eqref{eq:forward-backward representation}
with $i=2$, one obtains
\[
    \cE(u)
    =
    x_u(T)
    =
    \Phi_{t_0,T}(x^0)
    +
    \int_{t_0}^T
    D\Phi_{t,T}(x_u(t))B(t,x_u(t))u(t)\,dt
    =
    G_2(u)+\Phi_{t_0,T}(x^0)+y_2,
\]
and therefore $D\cE(u)=DG_2(u)$. The case $i=1$ follows similarly from
the same representation, which gives $DG_{\tau_i}(u)
    =
    D\Phi_{T,\tau_i}(x_u(T))\,D\cE(u)$.
This proves~\eqref{eq:identities}.

Next, suppose that $\cN_{\tau_i}(u)$ is invertible. Then, the canonical right inverse of $L_{u,\tau_i}$ satisfies
\[
R_{u,\tau_i}:=L_{u,\tau_i}^\ast\cN_{\tau_i}(u)^{-1}\in\mathscr{L}(\R^d,\cX),\qquad\,L_{u,\tau_i}R_{u,\tau_i}=\idty.
\]
Hence, if $\idty+A_{u,\tau_i}$ is invertible with $A_{u,\tau_i}:=K_{u,\tau_i}L_{u,\tau_i}^\ast\cN_{\tau_i}(u)^{-1}$, one uses $DG_{\tau_i}(u)=L_{u,\tau_i}+K_{u,\tau_i}$ to get
\begin{equation}\label{eq:idenditie for right inverse}
    (L_{u,\tau_i}+K_{u,\tau_i})\,R_{u,\tau_i}\,(\idty+A_{u,\tau_i})^{-1}=\idty\quad\text{which is equivalent to}\quad DG_{\tau_i}(u)\,L_{u,\tau_i}^\ast\begin{bmatrix}
    DG_{\tau_i}(u)\,L_{u,\tau_i}^\ast
\end{bmatrix}^{-1}=\idty
\end{equation}
showing that $DG_{\tau_i}(u)\,L_{u,\tau_i}^\ast\in\R^{d\times d}$ is invertible and $DG_{\tau_i}(u)$ is right-invertible. So $DG_{\tau_i}(u)$ is onto, which implies that $M(u):=DG_{\tau_i}(u)DG_{\tau_i}(u)^\ast$ is invertible.

Suppose now that $M(u):=DG_{\tau_i}(u)DG_{\tau_i}(u)^\ast$ is invertible. The canonical right-inverse of $DG_{\tau_i}(u)$ satisfies
\[
T_{u,\tau_i}:=DG_{\tau_i}(u)^\ast M(u)^{-1}\in\mathscr{L}(\R^d,\cX),\qquad\,DG_{\tau_i}(u)T_{u,\tau_i}=\idty.
\]
If $\idty-A_{u,\tau_i}$ is invertible with $A_{u,\tau_i}:=K_{u,\tau_i}DG_{\tau_i}(u)^\ast M(u)^{-1}$, one uses
$L_{u,\tau_i}=DG_{\tau_i}(u)-K_{u,\tau_i}$ to obtain
\begin{equation}\label{eq:important}
    (DG_{\tau_i}(u)-K_{u,\tau_i})\,T_{u,\tau_i}\,(\idty-A_{u,\tau_i})^{-1}=\idty\quad\Longleftrightarrow\quad L_{u,\tau_i}\,DG_{\tau_i}(u)^\ast\begin{bmatrix}
    L_{u,\tau_i}\,DG_{\tau_i}(u)^\ast
\end{bmatrix}^{-1}=\idty.
\end{equation}
This shows that $L_{u,\tau_i}\,DG_{\tau_i}(u)^\ast\in\R^{d\times d}$ is invertible and $L_{u,\tau_i}$ is right-invertible. So $L_{u,\tau_i}$ is onto, which implies that $\cN_{\tau_i}(u):=L_{u,\tau_i}\,L_{u,\tau_i}^\ast$ is invertible. This completes the proof of the lemma.

\subsection{Proof of Theorem~\ref{thm:minimizer-is-fp}}\label{ss:minimizer-is-fp}
In this section, we use the same notations introduced in Section~\ref{ss:preparatory result}.
\begin{proof}[\textit{Proof} of  Theorem~\ref{thm:minimizer-is-fp}]
Since $\idty+K_{u,\tau_i}L_{u,\tau_i}^\ast\cN_{\tau_i}(u)^{-1}$ and $\cN_{\tau_i}(u)$ are invertible for $u\in\cF(y_i)\subset\cF(C_i)\subset L^\infty\subset L^2$, one deduces that $DG_{\tau_i}(u)$ is onto by Lemma~\ref{lem:preparatory result}. Then, for any local minimizer $\bar u$ on $\mathfrak F_i:=\{u\in\cF(y_i):\ G_{\tau_i}(u)=0\}$, the Lagrange multiplier theorem in Hilbert spaces~\cite[Theorem 43.D, p.~290]{zeidler2013nonlinear}
(applied with $F(u)=\tfrac12\|u\|_2^2$ and the submersion $G_{\tau_i}:\cF(y_i)\to\R^d$, $G_{\tau_i}(u)=L_{u,\tau_i}u-y_i$) yields $\lambda\in\R^d$ such that
\begin{equation}\label{eq:u bar}
    \bar u=DG_{\tau_i}(\bar u)^\ast\lambda,\qquad L_{\bar u,\tau_i}\bar u=y_i.
\end{equation}
Since $DG_{\tau_i}(\bar u)$ is onto, one has that $M(\bar u)=DG_{\tau_i}(\bar u)DG_{\tau_i}(\bar u)^\ast\in\R^{d\times d}$ is invertible, and Lemma~\ref{lem:preparatory result} ensures that $L_{\bar u,\tau_i}\,DG_{\tau_i}(\bar u)^\ast\in\R^{d\times d}$ is invertible. Next, left-multiplying the first identity in~\eqref{eq:u bar} by $L_{\bar u,\tau_i}$ and using the second one yields
\begin{equation}
    \lambda = \begin{bmatrix}
        L_{\bar u,\tau_i}\,DG_{\tau_i}(\bar u)^\ast
    \end{bmatrix}^{-1}y_i
\end{equation}
Substituting back to~\eqref{eq:u bar} leads to
\begin{equation}\label{eq:u bar 1}
    \bar u = P_{\bar u,\tau_i}y_i\qquad\text{where}\quad P_{\bar u,\tau_i}:=DG_{\tau_i}(\bar u)^\ast\begin{bmatrix}
        L_{\bar u,\tau_i}\,DG_{\tau_i}(\bar u)^\ast
    \end{bmatrix}^{-1}.
\end{equation}
Recall from~\eqref{eq:important} that 
\[
L_{\bar u,\tau_i}\,P_{\bar u,\tau_i}=\idty
\]
showing that $P_{\bar u,\tau_i}$ is a right-inverse of $L_{\bar u,\tau_i}$.
Next, since the canonical right inverse $R_{\bar u,\tau_i}=L_{\bar u,\tau_i}^{\ast}\cN_{\tau_i}(\bar u)^{-1}$ of $L_{\bar u,\tau_i}$ satisfies $L_{\bar u,\tau_i}R_{\bar u,\tau_i}=\idty$, one deduces that 
$L_{\bar u,\tau_i}\,(P_{\bar u,\tau_i}-R_{\bar u,\tau_i})=0$, which implies
\[
(P_{\bar u,\tau_i}-R_{\bar u,\tau_i})y\in\ker L_{\bar u,\tau_i}\qquad\forall y\in\R^d.
\]
Therefore, \eqref{eq:u bar 1} recast as
\begin{equation}\label{eq:u bar decomposition}
    \bar u =  R_{\bar u,\tau_i}y_i + z_i,\qquad z_i:=(P_{\bar u,\tau_i}-R_{\bar u,\tau_i})y_i\in\ker L_{\bar u,\tau_i}.
\end{equation}

Now, let us prove the equivalence in the theorem. Assume that
\[
\langle[L_{\bar u,\tau_i}DG_{\tau_i}(\bar u)^\ast]^{-1}y_i,DG_{\tau_i}(\bar u)\,h\rangle_{\R^d}=0\qquad\forall\,h\in\ker L_{\bar u,\tau_i}
\]
On one hand, for $h\in\ker\,L_{\bar u,\tau_i}$, one get from~\eqref{eq:u bar 1} that $\langle\bar u,h\rangle_{L^2} = \langle[L_{\bar u,\tau_i}DG_{\tau_i}(\bar u)^\ast]^{-1}y_i,DG_{\tau_i}(\bar u)\,h\rangle_{\R^d} = 0$, which 
implies that
\begin{equation}\label{eq:prop of bar u}
    \bar u\in(\ker\,L_{\bar u,\tau_i})^\perp.
\end{equation} 
On the other hand, by~\eqref{eq:u bar decomposition}, \eqref{eq:prop of bar u} and $R_{\bar u,\tau_i}y_i=L_{\bar u,\tau_i}^{\ast}\cN_{\tau_i}(\bar u)^{-1}y_i\in\Ran(L_{\bar u,\tau_i}^{\ast})=(\ker L_{\bar u,\tau_i})^\perp$, one finds that $z_i=\bar u-R_{\bar u,\tau_i}y_i\in (\ker L_{\bar u,\tau_i})^\perp$. Therefore, $z_i\in \ker(L_{\bar u,\tau_i})\cap(\ker L_{\bar u,\tau_i})^\perp=\{0\}$ so $z_i=0$, and 
\[
\bar u=R_{\bar u,\tau_i}y_i=L_{\bar u,\tau_i}^{\ast}\cN_{\tau_i}(\bar u)^{-1}y_i=\cS_i(\bar u).
\]
Conversely, if $\bar u=L_{\bar u,\tau_i}^{\ast}\cN_{\tau_i}(\bar u)^{-1}y_i$, then~\eqref{eq:u bar decomposition} yields $P_{\bar u,\tau_i}y_i=R_{\bar u,\tau_i}y_i\in\Ran(L_{\bar u,\tau_i}^{\ast})=(\ker L_{\bar u,\tau_i})^\perp$. Let $h\in\ker\,L_{\bar u,\tau_i}$, then 
\[
\langle[L_{\bar u,\tau_i}DG_{\tau_i}(\bar u)^\ast]^{-1}y_i,DG_{\tau_i}(\bar u)\,h\rangle_{\R^d}=\langle[DG_{\tau_i}(\bar u)^\ast L_{\bar u,\tau_i}DG_{\tau_i}(\bar u)^\ast]^{-1}y_i\,h\rangle_{L^2}=\langle P_{\bar u,\tau_i}y_i,h\rangle_{L^2}=0.
\]
This completes the proof of the theorem.
\end{proof}

\subsection{Proof of Proposition~\ref{prop:exist-on-ball-global}}\label{ss:exist-on-ball-global}
\begin{proof}
For ease in notation, set $L^p:=L^p((t_0, T); \R^k)$ for $p\in\{1,2,\infty\}$. By definition, $\mathfrak F_i\subset\,L^\infty\subset\,L^2$. Let $u\in\,L^2$ and $(u_n)\subset\mathfrak F_i$ are such that $u_n\rightharpoonup u$ in $L^2$. Let us show that $u\in\mathfrak F_i$.

\paragraph{\textbf{Step 1}} One has $u\in\,L^\infty$ and $\|u\|_\infty\le\,\zeta_i$. In fact, since $(u_n)\subset\,L^\infty$, by Banach--Alaoglu, there exists a subsequence (not relabeled) and $\tilde u\in L^\infty$ with 
$\|\tilde u\|_{L^\infty}\le \zeta_i$ such that $u_n\stackrel{\ast}{\rightharpoonup}\tilde u$ in $L^\infty$, i.e.
$\int_{t_0}^T \varphi(t)^\top\, u_n(t)\,dt\to\int_{t_0}^T \varphi(t)^\top\tilde u(t)\,dt$ for all $\varphi\in L^1$.
On the other hand, $u_n\rightharpoonup u$ in $L^2$, so $\int_{t_0}^T \psi(t)^\top u_n(t)\,dt\to\int_{t_0}^T \psi(t)^\top u(t)\,dt$ for all $\psi\in L^2$. It holds that
\[
\int_{t_0}^T \phi(t)^\top u(t)\,dt=\lim_{n\to\infty}\int_{t_0}^T \phi(t)^\top u_n(t)\,dt=\int_{t_0}^T \phi(t)^\top \tilde u(t)\,dt\qquad\forall\phi\in\,L^2\cap\,L^1=L^2.
\]
In particular, for $\phi:=u-\tilde{u}\in\,L^2$, one finds $\int_{t_0}^T|u(t)-\tilde{u}(t)|^2\,dt=0$, which implies $u=\tilde{u}$, a.e., and thus $u\in L^\infty$ and $\|u\|_{L^\infty}=\|\tilde u\|_{L^\infty}\le \zeta_i$.

\paragraph{\textbf{Step 2}} One has $G_{\tau_i}(u):=L_{u,\tau_i}u-y_i=0$. We split the proof of this fact into three steps.

\paragraph{\textbf{Step 2.1}} For $n\in\N$, and $t\in[t_0, T]$, set $\delta_n(t)=\int_{t_0}^tB(s,x_u(s))\,(u_n(s)-u(s))\,ds$. It is clear that $(\delta_n)\subset\,C^0([t_0, T]; \R^d)$. Consider the bounded and linear operator $T_t:L^2\to\R^d$, $T_th=\int_{t_0}^T\mathbf 1_{[t_0,t]}(s)\,B(s,x_u(s))\,h(s)\,ds$. Its adjoint $T_t^\ast:\R^d\to\,L^2$ is $T_t^\ast\,z = \mathbf 1_{[t_0,t]}(\cdot)\,B(\cdot,x_u(\cdot))^\top\,z\in\,L^2$ for all $z\in\R^d$. Then $\delta_n(t)=T_t(u_n-u)$ and the weak convergence in $L^2$ gives $\langle\delta_n(t),z\rangle_{\R^d}=\langle\,u_n-u,T_t^\ast\,z\rangle_{L^2}\to\,0$ for all $z\in\R^d$, hence $\delta_n(t)\to\,0$ in $\R^d$.

 For all $n\in\N$, one has by Cauchy-Schwarz inequality,
\[
|\delta_n(t)|\le\|B\|_\infty\,(\zeta_i+\|u\|_2)\,\sqrt{\Delta t_0}\qquad\forall\,t\in[t_0, T]
\]
showing that $(\delta_n)_{n}$ is uniformly bounded. For $(t, t')\in[t_0, T]^2$ with $t\ge t'$ (the case of $t<t'$ is identical), one has by Cauchy-Schwarz inequality,
\[
|\delta_n(t)-\delta_n(t')|\le\|B\|_\infty\,(\zeta_i+\|u\|_2)\,\sqrt{t-t'}\qquad\forall\,n\in\N
\]
showing that $(\delta_n)_{n}$ is equicontinuous. Consequently, $(\delta_n)\subset\,C^0([t_0, T]; \R^d)$ is relatively compact, and by Ascoli-Arzèla, it admits a uniformly convergent subsubsequence; its limit must be the pointwise limit, viz.
\[
\|\delta_n\|_\infty\to\,0\quad\text{as}\quad\,n\to\infty.
\]

\paragraph{\textbf{Step 2.2}} Let $x_{u_n}$ and $x_u$ be the solution of $\dot x = N_t(x)+B(t,x)\,u_n$ and $\dot x = N_t(x)+B(t,x)\,u$, respectively with the same intitial state $x^0\in\R^d$. Gronwall's lemma yields
\[
\|x_{u_n}-x_u\|_\infty\le\Delta t_0\,e^{(\Lambda_1+L_B\zeta_i)\Delta t_0}\|\delta_n\|_\infty\to\,0\quad\text{as}\quad\,n\to\infty.
\]

\paragraph{\textbf{Step 2.3}}  Write $K_i(u)(t):=D\Phi_{t,\tau_i}(x_u(t))\,B(t,x_u(t))\in\R^{d\times k}$. Lemma~\ref{lem:spectral norm of D_s_t_phi and D^2_s_t_phi} and~\eqref{eq:cst1} yield 
\[
\|K_i(u_n)-K_i(u)\|_\infty\le\,\bigg(\frac{\Lambda_2}{\Lambda_1}\,E_{2,\infty}\|B\|_\infty+L_B\,E_{1,\infty}\bigg)\|x_{u_n}-x_u\|_\infty\to\,0\quad\text{as}\quad\,n\to\infty.
\]
Overall, 
\[
\begin{aligned}
G_{\tau_i}(u_n)-G_{\tau_i}(u)
&=\int_{t_0}^T\!\big(K_i(u_n)-K_i(u)\big)\,u_n\,dt\;+\;\int_{t_0}^T\!K_i(u)\,(u_n-u)\,dt \\
&\xrightarrow[n\to\infty]{} 0+0,
\end{aligned}
\]
since the first term vanishes by $\|K_i(u_n)-K_i(u)\|_{L^\infty}\to0$ and $L^2$-boundedness of $(u_n)$, and the second by weak convergence against the fixed $K_i(u)\in L^\infty\subset\,L^2$. As $G_{\tau_i}(u_n)=0$, we obtain $G_{\tau_i}(u)=0$.

\paragraph{\textbf{Step 3}} One has $\lambda_{\min}\big(\cN_{\tau_i}(u)\big)\ge C_i^{-1}$. In fact, $\cN_{\tau_i}(u)=\int_{t_0}^T K_i(u)(t)\,K_i(u)(t)^\top\,dt$, and spectral stability,
\[
\big|\lambda_{\min}(\cN_{\tau_i}(u_n))-\lambda_{\min}(\cN_{\tau_i}(u))\big|
\ \le\|\cN_{\tau_i}(u_n)-\cN_{\tau_i}(u)\|\le 2\Delta t_0\,\|B\|_\infty\,e^{\Lambda_1\Delta t_0}\,\|K_i(u_n)-K_i(u)\|_\infty \ \xrightarrow[n\to\infty]{} 0.
\]
Since $\lambda_{\min}(\cN_{\tau_i}(u_n))\ge C_i^{-1}$ for all $n$, we obtain
$\lambda_{\min}(\cN_{\tau_i}(u))\ge C_i^{-1}$.

If $(u_n)\subset\mathfrak F_i$ is a minimizing sequence, then $(\tfrac12\|u_n\|_{L^2}^2)_n$ is a bounded sequence. After extraction, $u_{n_j}\rightharpoonup \bar u$ in $L^2$.
By weak sequential closedness, $\bar u\in\mathfrak F_i$.
Then $\tfrac12\|\bar u\|_{L^2}^2 \le \liminf_j \tfrac12\|u_{n_j}\|_{L^2}^2$,
so $\bar u$ attains the infimum. 
\end{proof}

\subsection{Proof of Proposition~\ref{pro:uniformly invertible B}}\label{ss:uniformly invertible B}

\begin{proof}
    Let $(x^0,\,u)\in\R^d\times L^\infty((t_0, T); \R^k)$ and $x_u$ be the solution to~\eqref{eq:state-dependent-input}. Then, applying~\eqref{eq:EGectral norm of D_s_t_phi} with $\beta(t)=\Phi_{t,t_0}(x_u(t))$, one immediately finds that
    \[
    |D\Phi_{t,\tau_i}(x_u(t))z|\ge e^{-\Lambda_1|\tau_i-t|}|z|\qquad\forall(t,\,z)\in[t_0, T]\times\R^d,\quad\tau_1=t_0,\,\tau_2=T.
    \]
    It follows that
    \begin{equation*}
        y^\top \cN_{\tau_i}(u)\,y=\int_{t_0}^T |D\Phi_{t,\tau_i}(x_u(t))^\top\,B(t,x_u(t))^\top y|^2\,dt\ge e^{-2\Lambda_1\Delta t_0}\,|y|^2\int_{t_0}^Tb(t)^2\,dt\ge \frac{e^{-2\Lambda_1\Delta t_0}\,\|b\|_1^2}{\Delta t_0}\,|y|^2,
    \end{equation*}
by the Cauchy-Schwarz inequality. This proves~\eqref{eq:coercivity bound full actuated case}. So, letting $C_i>0$ as in the proposition, one gets
\[
\cF(C_i) = \left\{u\in L^\infty((t_0, T); \R^k):\lambda_{\min}\big(\cN_{\tau_i}(u)\big)\ge\,e^{-2\Lambda_1\Delta t_0}\|b\|_1^2/\Delta t_0\right\}=L^\infty((t_0, T); \R^k),
\]
and \eqref{eq:self-mapping} is automatically satisfied.
\end{proof}

\subsection{Proof of Theorem~\ref{thm:relative coercivity from a control reference}}\label{ss:relative coercivity from a control reference}

\begin{proof}
    Let $C_i>0$ as in the statement of the proposition, and $y_i\in\R^d$. Then, \eqref{eq:self-mapping} is satisfied if and only if $v\in\cF(C_i)$ for every $v\in\cS_i(\cF(y_i))$. Observe that
\[
\cF(C_i) = \left\{u\in L^\infty((t_0, T); \R^k):\lambda_{\min}\big(\cN_{\tau_i}(u)\big)\ge\lambda_{\min}\big(\cN_{\tau_i}(u_i^\flat)\big)/(1+\theta)\right\}.
\]
Now, $v\in\cS_i(\cF(y_i))$ if and only if $v=\cS_i(u)=B(t,x_u(t))^\top D\Phi_{t,\tau_i}(x_{u}(t))^\top\,\cN_{\tau_i}(u)^{-1}\,y_i$ for some $u\in\cF(C_i)$, where $\tau_1:=t_0$ and $\tau_2:=T$. Therefore, $\|v\|_\infty\le\|B\|_\infty\,e^{\Lambda_1\Delta t_0}\lambda_{\min}^{-1}(\cN_{\tau_i}(u))\,|y_i|\le(1+\theta)\,e^{\Lambda_1\Delta t_0}\|B\|_\infty|y_i|/\lambda_{\min}\big(\cN_{\tau_i}(u_i^\flat)\big)$.
It remains to prove that 
\begin{equation}\label{eq:to show}
\lambda_{\min}(\cN_{\tau_i}(v))=\lambda_{\min}\big(\cN_{\tau_i}(\cS_i(u))\big)\ge\lambda_{\min}\big(\cN_{\tau_i}(u_i^\flat)\big)/(1+\theta).     
\end{equation}
By the spectral stability and $L_{\cN_{\tau_i}}$-Lipschitz of $u\in\{u\in\,L^\infty((t_0, T); \R^k):\|u\|\le\,\zeta_i \}\mapsto\cN_{\tau_i}(u)$, one finds that
\begin{equation*}
     \left|\lambda_{\min}\big(\cN_{\tau_i}(\cS_i(u))\big)-\lambda_{\min}\big(\cN_{\tau_i}(u_i^\flat)\big)\right|\le L_{\cN_{\tau_i}}\left(\frac{(1+\theta)\,e^{\Lambda_1\Delta t_0}\|B\|_\infty|y_i|}{\lambda_{\min}\big(\cN_{\tau_i}(u_i^\flat)\big)}+\|u_i^\flat\|_\infty\right).
\end{equation*}
Since 
\[
\left|\lambda_{\min}\big(\cN_{\tau_i}(\cS_i(u))\big)-\lambda_{\min}\big(\cN_{\tau_i}(u_i^\flat)\big)\right|\ge\frac{\lambda_{\min}\big(\cN_{\tau_i}(u_i^\flat)\big)}{1+\theta}+\frac{\theta\lambda_{\min}\big(\cN_{\tau_i}(u_i^\flat)\big)}{1+\theta}-\lambda_{\min}\big(\cN_{\tau_i}(\cS_i(u))\big),
\]
one deduces that if 
\begin{equation}\label{eq:sufficient conditions for invertibility general proof}
    \frac{\theta\lambda_{\min}\big(\cN_{\tau_i}(u_i^\flat)\big)}{1+\theta}\ge L_{\cN_{\tau_i}}\left(\frac{(1+\theta)\,e^{\Lambda_1\Delta t_0}\|B\|_\infty|y_i|}{\lambda_{\min}\big(\cN_{\tau_i}(u_i^\flat)\big)}+\|u_i^\flat\|_\infty\right)
\end{equation}
then~\eqref{eq:to show} is satisfied. Since~\eqref{eq:sufficient conditions for invertibility general proof} is equivalent to~\eqref{eq:sufficient conditions for invertibility general}, this completes the proof of the proposition.
\end{proof}

\subsection{Proof of Lemma~\ref{lem:Lipschitz constant of N_i}}\label{ss:proof of Lipschitz constant of N_i}

\begin{proof}
    Write $K_i(u)(t):=D\Phi_{t,\tau_i}(x_u(t))\,B(t,x_u(t))\in\R^{d\times k}$ so that $\cN_{\tau_i}(u)= \int_{t_0}^T K_i(u)(t)\,K_i(u)(t)^\top\,dt$.
    Clearly, $\cN_{\tau_i}$ is continuous on $L^\infty$ by Lebesgue dominated convergence theorem. Let us prove the second part of the lemma. Let $u,\,v\in\{u\in\,L^\infty((t_0, T); \R^k):\|u\|_\infty\le\,\zeta_i \}$. By the Cauchy-Schwarz inequality, one finds that
    \[
    \frac{d}{dt}|x_u(t)-x_v(t)|=|\dot{x}_u(t)-\dot{x}_v(t)|\le(\Lambda_1+L_B\,\zeta_i)|x_u(t)-x_v(t)|+\|B\|_\infty\,|u(t)-v(t)|
    \]
    which, by Gronwall's lemma, yields
    \begin{equation}\label{eq:estimate on |x_u(t)-x_v(t)|}
        |x_u(t)-x_v(t)|\le\|B\|_\infty\frac{e^{(\Lambda_1+L_B\,\zeta_i)(t-t_0)}-1}{(\Lambda_1+L_B\,\zeta_i)}\|u-v\|_\infty\qquad\forall \in[t_0, T].
    \end{equation}
    Using~Lemma~\ref{lem:spectral norm of D_s_t_phi and D^2_s_t_phi}, one finds (we let $\Delta t_0:=T-t_0$)
    \begin{eqnarray*}\label{eq:inter N_2(u)-N_2(v)}
       \frac{\|\cN_{\tau_2}(u)-\cN_{\tau_2}(v)\|}{\|u-v\|_\infty}
       &\le&\frac{2\|B\|_\infty^3\Lambda_2}{\Lambda_1(\Lambda_1+L_B\,\zeta_2)}\int_{t_0}^T\,e^{\Lambda_1(T-t)}\,\left(e^{2\Lambda_1(T-t)}-e^{\Lambda_1(T-t)}\right)\,(e^{(\Lambda_1+L_B\,\zeta_2)(t-t_0)}-1)\,dt\nonumber\\
       & &+\frac{2L_B\|B\|_\infty^2}{(\Lambda_1+L_B\,\zeta_2)}\int_{t_0}^Te^{2\Lambda_1(T-t)}(e^{(\Lambda_1+L_B\,\zeta_2)(t-t_0)}-1)\,dt\nonumber\\
       &=&\frac{L_B\|B\|_\infty^2e^{\Lambda_1\Delta t_0}}{\Lambda_1}\left[\frac{e^{\Lambda_1\Delta t_0}-e^{L_B\zeta_2\Delta t_0}}{\Lambda_1-L_B\zeta_i}-\frac{e^{L_B\zeta_2\Delta t_0}-e^{-\Lambda_1\Delta t_0}}{L_B\zeta_2+\Lambda_1}\right]+\nonumber\\
       &&\hspace{-2cm}\frac{2\|B\|_\infty^3\Lambda_2e^{\Lambda_1\Delta t_0}}{\Lambda_1(\Lambda_1+L_B\,\zeta_2)}\left[\frac{e^{2\Lambda_1\Delta t_0}-e^{L_B\zeta_2\Delta t_0}}{2\Lambda_1-L_B\zeta_2}+\frac{e^{\Lambda_1\Delta t_0}-e^{L_B\zeta_2\Delta t_0}}{\Lambda_1-L_B\zeta_2}+\frac{e^{\Lambda_1\Delta t_0}-e^{-\Lambda_1\Delta t_0}}{2\Lambda_1}-\frac{e^{2\Lambda_1\Delta t_0}-e^{-\Lambda_1\Delta t_0}}{3\Lambda_1}\right],
    \end{eqnarray*}
    by a direct integration. The same proof applies to $\cN_{\tau_1}$.
\end{proof}

\subsection{Proof of Proposition~\ref{pro:Lipschitz constant of N_i for HNN}}\label{ss:proof of Lipschitz constant of N_i for HNN}

\begin{proof}
    Recall that $DN(x)=-D+W\,D\sigma(x)$ 
    and $D^2N(x) = W\,D^2\sigma(x)$.
    Then, from
    \[
    \frac{d}{dt}\phi_t(x) = N(\phi_t(x)),\quad\phi_0(x)=x;\qquad\frac{d}{dt}D\phi_t(x) = DN(\phi_t(x))\,D\phi_t(x),\quad D\phi_0(x)=\idty\quad\forall x\in\R^d,
    \]
    Cauchy-Schwarz inequality yields
    \[
     \frac{d}{dt}|D\phi_t(x)\,y| = \frac{\langle\,-D\,D\phi_t(x)\,y,D\phi_t(x)\,y\rangle+\langle\,W\,D\sigma(\phi_t(x))\,D\phi_t(x)\,y,D\phi_t(x)\,y\rangle}{|D\phi_t(x)\,y|}\le\Gamma\,|D\phi_t(x)\,y|\quad\forall y\in\R^d\bs\{0\},
    \]
    which, by Gronwall's lemma, leads to
    \begin{equation}\label{eq:for Dflow}
        \|D\phi_t(x)\|\le e^{\Gamma t}\qquad\forall\,x\in\R^d,\,\forall t\ge 0.
    \end{equation}
    Similarly, one has
    \[
    \frac{d}{dt}D^2\phi_t(x) =  D^2N(\phi_t(x))\,\left[D\phi_t(x),D\phi_t(x)\right]+DN(\phi_t(x))\,D^2\phi_t(x),\quad D^2\phi_0(x) = 0,\quad x\in\R^d,
    \]
    so that letting $g(t):=D^2\phi_t(x)[y,z]$ for $y, z\in\R^d$, Cauchy-Schwarz inequality and~\eqref{eq:for Dflow} lead to
    \begin{eqnarray*}
          \frac{d}{dt}|g(t)|&=&\frac{\langle\,-D\,g(t),g(t)\rangle+\langle\,W\,D\sigma(\phi_t(x))\,g(t),g(t)\rangle+\langle\,W\,D^2\sigma(\phi_t(x))\left[D\phi_t(x)y,D\phi_t(x)z\right],g(t)\rangle}{|g(t)|}\nonumber\\
          &\le&\Gamma\,|g(t)|+\Gamma_2\,e^{2\Gamma\,t}\,|y|\,|z|,
    \end{eqnarray*}
    which, by Gronwall's lemma and immediate integration, yields
    \begin{equation}\label{eq:for D2flow}
        \|D^2\phi_t(x)\|\le\begin{cases}
        \Gamma_2\,\frac{e^{2\Gamma\,t}-e^{\Gamma\,t}}{\Gamma}&\quad\text{if}\quad\Gamma\neq 0,\cr
        \Gamma_2\,t&\quad\text{if}\quad\Gamma= 0,
    \end{cases}\qquad\forall\,x\in\R^d,\,\forall t\ge 0.
    \end{equation}
    In particular as in Lemma~\ref{lem:spectral norm of D_s_t_phi and D^2_s_t_phi}, one can replace $x\in\R^d$ in~\eqref{eq:for Dflow} and~\eqref{eq:for D2flow} with $\beta(\cdot)$ for any $\beta\in C^0([0, T];\R^d)$.  Finally, for $u,v\in L^\infty((0,T);\R^k)$, one uses the same technique to prove~\eqref{eq:for Dflow} and~\eqref{eq:for D2flow}, and get
    \begin{equation}\label{eq:|x_u-x_v| HNN}
        |x_u(t)-x_v(t)|\le\begin{cases}
        \|B\|_\infty\,\frac{e^{\Gamma\,t}-1}{\Gamma}\|u-v\|_\infty&\quad\text{if}\quad\Gamma\neq 0,\cr
        \|B\|_\infty\,t\,\|u-v\|_\infty&\quad\text{if}\quad\Gamma= 0,
    \end{cases}\qquad\forall t\ge 0.
    \end{equation}

    Now, recall that $\cN_{\tau_2}(u)= \int_0^T\,D\phi_{t}(x_u(T-t))\,B(T-t)\,B(T-t)^\top\,D\phi_{t}(x_u(T-t))^\top\,dt$, so that
    letting $Q_u(t):=D\phi_{t}(x_u(T-t))$, if $\Gamma=0$, the estimates~\eqref{eq:for Dflow}, \eqref{eq:for D2flow} and~\eqref{eq:|x_u-x_v| HNN} yield
    \begin{equation}
         \|\cN_{\tau_2}(u)-\cN_{\tau_2}(v)\|\le\,2\,\|B\|_\infty^3\,\Gamma_2\|u-v\|_\infty\,\int_0^T\,t\,(T-t)\,dt = \frac{\Gamma_2\,\|B\|_\infty^3\,T^3}{3}\,\|u-v\|_\infty.
    \end{equation}
    When $\Gamma\neq 0$, one finds that 
    \begin{equation*}
         \|\cN_{\tau_2}(u)-\cN_{\tau_2}(v)\|\le\,2\,\|B\|_\infty^3\,\Gamma_2\,\|u-v\|_\infty\,\int_0^T\,\frac{(e^{3\Gamma\,t}-e^{2\Gamma\,t})(e^{\Gamma(T-t)}-1)}{\Gamma^2}\,dt=\frac{\Gamma_2\,\|B\|_\infty^3}{3}\left(\frac{e^{\Gamma\,T}-1}{\Gamma}\right)^3\,\|u-v\|_\infty,
    \end{equation*}
    which completes the proof of the inequality involving $\cN_{\tau_2}$. The same arguments are used to prove those of $\cN_{\tau_1}$.
\end{proof}

\bibliographystyle{siamplain}
\bibliography{references}
\end{document}